\documentclass{article}
\usepackage{arxiv}
\usepackage[utf8]{inputenc}
\usepackage[T1]{fontenc}
\usepackage{hyperref}
\usepackage{url}
\usepackage{booktabs}
\usepackage{amsfonts}
\usepackage{nicefrac}
\usepackage{microtype}
\usepackage{lipsum}
\usepackage{graphicx}
\usepackage{dcolumn}
\usepackage{bm}
\usepackage{amsmath}
\usepackage{amssymb}
\usepackage{amscd}
\usepackage{amsthm}
\usepackage{amsfonts}
\usepackage{graphicx}
\usepackage{multirow}

\newcommand{\beq}{\begin{equation}}
\newcommand{\eeq}{\end{equation}}
\newcommand{\bea}{\begin{eqnarray}}
\newcommand{\eea}{\end{eqnarray}}
\newcommand{\bc}{\begin{cases}}
\newcommand{\ec}{\end{cases}}

\newtheorem{theorem}{Theorem}

\newtheorem{corollary}[theorem]{Corollary}

\newtheorem{definition}[theorem]{Definition}

\newtheorem{lemma}[theorem]{Lemma}

\begin{document}

\date{\today}
\author{ Jos\'{e} M. Amig\'{o} \\
Centro de Investigaci\'{o}n Operativa, \\
Universidad Miguel Hern\'{a}ndez, \\
03202 Elche, Spain \\
\texttt{jm.amigo@umh.es} \And Karsten Keller \\
Institut f\"{u}r Mathematik, \\
Universit\"{a}t zu L\"{u}beck, \\
23562 L\"{u}beck, Germany \\
\texttt{unakafova@math.uni-luebeck.de} \And Valentina Unakafova \\
Graduate School for Computing in Medicine and Life Science, \\
Universit\"{a}t zu L\"{u}beck, \\
23562 L\"{u}beck, Germany \\
\texttt{unakafova@math.uni-luebeck.de} }
\title{On entropy, entropy-like quantities, and applications}
\maketitle

\begin{abstract}
This is a review on entropy in various fields of mathematics and science.
Its scope is to convey a unified vision of the classical as well as some
newer entropy notions to a broad audience with an intermediate background in
dynamical systems and ergodic theory. Due to the breadth and depth of the
subject, we have opted for a compact exposition whose contents are a
compromise between conceptual import and instrumental relevance. The
intended technical level and the space limitation born furthermore upon the
final selection of the topics, which cover the three items named in the
title. Specifically, the first part is devoted to the avatars of entropy in
the traditional contexts: many particle physics, information theory, and
dynamical systems. This chronological order helps to present the materials
in a didactic manner. The axiomatic approach will be also considered at this
stage to show that, quite remarkably, the essence of entropy can be
encapsulated in a few basic properties. Inspired by the classical entropies,
further akin quantities have been proposed in the course of time, mostly
aimed at specific needs. A common denominator of those addressed in the
second part of this review is their major impact on research. The final part
shows that, along with its profound role in the theory, entropy has
interesting practical applications beyond information theory and
communications technology. For this sake we preferred examples from applied
mathematics, although there are certainly nice applications in, say,
physics, computer science and even social sciences. This review concludes
with a representative list of references.
\end{abstract}

\tableofcontents


\section{Introduction}

\label{sec:1}

Entropy is a general concept that appears in different settings with
different meanings. Thus, the Boltzmann-Gibbs entropy measures the
microscopic disorder in statistical mechanics; the Shannon entropy measures
uncertainty and compression in information theory; the Kolmogorov-Sinai
entropy measures (pseudo-)randomness in measure-preserving dynamical
systems; and the topological entropy measures complexity in topological
dynamics. As for its importance, entropy is the protagonist of the second
law of thermodynamics, associated by some authors with the arrow of time. In
information theory, coding theory, and cryptography, Shannon entropy lies at
the core of the fundamental definitions (information, typicality, channel
capacity,...) and results (asymptotic equipartition, channel coding theorem,
channel capacity theorem, conditions for secure ciphers,...). And in ergodic
theory, Kolmogorov-Sinai and topological entropy are perhaps the most
important invariants of metric and topological conjugacy, respectively,
which are the equivalence concepts in measure-preserving and topological
dynamics.

Also remarkable is a sort of universality that entropy possesses. By this we
mean that other indicators of, say, compressibility or complexity turn out
to be related to, or even coincide with some of the standard entropies. This
is what happens, for example, with algorithmic complexity in computer
science, Lempel-Ziv complexity in information theory, finitely observable
invariants for the class of all finitely-valued ergodic processes in the
theory of stochastic processes (see Theorem \ref{Thm3.4} below), and
permutation entropy in dynamical systems. Shannon himself characterized his
function as the only one that satisfies three basic properties or axioms,
which would explain its uniqueness \cite{Shannon1948}. Later on, the
`entropies' introduced by R\'{e}nyi \cite{Renyi1961}, Havrda-Charv\'{a}t 
\cite{HavrdaCharvat1967}, and Tsallis \cite{Tsallis1988} drew the
researchers' attention to other solutions under weaker conditions or
different sets of axioms. Alone the suppression of the fourth
Shannon-Khinchin axiom (\textquotedblleft the entropy of a system ---split
into subsystems $A$ and $B$--- equals the entropy of $A$ plus the
expectation value of the entropy of $B$, conditional on $A$%
\textquotedblright ) opens the door to a broad variety of entropy-like
quantities \cite{Hanel2011}, called generalized entropies, which include the
aforementioned R\'{e}nyi, Havrda-Charv\'{a}t and Tsallis entropies.

History shows that statistical mechanics, information theory, and dynamical
systems build a circle of ideas. First of all, the theory of dynamical
systems is an outgrowth of statistical mechanics, from which it borrowed the
fundamental concepts: state space (read phase space), orbit (read
trajectory), invariant measure (read Liouville measure), ergodicity (read
Boltzmann's \textit{Ergodensatz}), stationary state (read equilibrium
state),... and entropy! Independently, Shannon created information theory in
a purely probabilistic setting. Indeed, information is a function of a
probability distribution and information sources are stationary random
processes. Yet, the solution of the famous Maxwell's paradox in statistical
mechanics (in which a microscopic `demon' perverts the order established by
the second law of thermodynamics) required the intervention of information
theory in the form of Landauer's principle \cite{Bennett2003}. Furthermore,
Shannon's brainchild inspired Kolmogorov's work on Bernoulli shifts. In
return, information theory has benefited from symbolic dynamics. At the
intersection of these three fields, amid past and on-going cross-pollination
among them, lie ergodic theory and entropy.

In the last decades new versions of entropy have come to the fore.
Approximate entropy \cite{EckmannRuelle1985,Pincus1991}, directional entropy 
\cite{Milnor1988}, ($\varepsilon,\tau$)-entropy \cite{Gaspard1993},
permutation entropy \cite{Bandt2002A}, sample entropy \cite{Takens1981},
transfer entropy \cite{Schreiber2000}, ... are some of the entropy-like
quantities proposed by researchers to cope with new challenges in time
series analysis, cellular automata, chaos, synchronization, multiscale
analysis, etc. Their relationship to the classical mathematical entropies
(i.e., the Shannon, Kolmogorov-Sinai, and topological entropies) is diverse.
Thus, some of them turn out to be equivalent to a classical counterpart
(e.g., metric and topological permutation entropy), or are defined with
their help (e.g., transfer entropy, a conditional mutual information). In
other cases, a new version can be considered a generalization of a classical
one (e.g., directional entropy, ($\varepsilon,\tau$)-entropy). Still in
other cases, it is a convenient approximation devised for a specific purpose
(e.g., approximate and sample entropy). We sidestep the question whether any
one of them is an `entropy' or rather an `entropy-like' quantity, and use
the word entropy in a wide sense. In any case, along with these new
proposals, the conventional entropies continue to play a fundamental role in
ergodic theory and applications. Always when it comes to capture the elusive
concepts of disorder, information (or ignorance), randomness, complexity,
irregularity, ..., whether in mathematics, natural sciences, or social
sciences, entropy has proved to be unrivaled.

In the remaining sections we are going to zoom in on most of the topics
mentioned above. At the same time we have tried to make our exposition
reasonably self-contained in regard to the basic concepts and results. For
details and more advanced materials the reader will be referred to the
literature. But beyond the formal issues, this review should be more than
just a guided visit to the most accessible facets of entropy. It should also
show that entropy, in whichever of its variants, is an exciting and useful
mathematical object.

\section{A panorama of entropy}

\label{section2}

This section is an overview of entropy in the most traditional settings. For
a more advanced account, which also includes interesting historical
information, the reader is referred to the excellent review \cite{Katok2007}
and the references therein.

\subsection{Entropy in thermodynamics and statistical mechanics}

The word \textit{entropy}\ was coined by the German physicist R. Clausius
(1822-1888) \cite{Clausius1865}, who introduced it in thermodynamics in 1865
to measure the amount of energy in a system that cannot produce work. The
fact that the entropy increases in irreversible adiabatic processes that
take one equilibrium state to another constitutes the second law of
thermodynamics and clearly shows the central role of entropy in the physics
of macroscopic bodies. The asymmetry in time introduced by this law has been
associated by some authors with the arrow of time.

Once the atomic nature of matter was discovered, physicists set out to find
the microscopic counterparts of all macroscopic quantities in general, and
entropy in particular. Think, for example, of one mole of an isolated gas in
equilibrium (e.g., one gram of hydrogen). Its macroscopic state or \textit{%
macrostate} is just determined by two parameters, say its volume $V$ and
temperature $T$. Microscopically the system consists of $N\simeq 6.022\times
10^{23}$ molecules moving (in a first approximation) under the laws of
classical point-mass mechanics. Alternatively one can think of the system as
a point describing a curve in a region $\Omega \subset \mathbb{R}^{6N}$,
called the \textit{phase space}, spanned by the coordinates and momenta of
the $N$ particles. It should be clear that the huge number of particles
makes an impracticable task to solve the equations of motion, while allowing
an statistical approach. Let $W$ be the number of microscopic states or 
\textit{microstates} (as determined by, e.g., the initial positions and
momenta of the particles) consistent with the imposed macroscopic
constraints ($V=V_{0}$, and $T=T_{0}$ in our example). For simplicity (and
in accordance with quantum physics!), we assume here that $W$ is finite. In
fact, the discretization of the state space was used by Max Planck as a
working hypothesis (following previous ideas of Boltzmann and before
Heisenberg's uncertainty principle justified it) to derive the radiation
spectrum of the black body, thus ushering in the new era of quantum
mechanics. The \textit{Boltzmann entropy} of the system is then%
\begin{equation}
S=k_{B}\ln W,  \label{Boltzmann_S}
\end{equation}%
where $k_{B}=1.3807\times 10^{-23}$ J/K is called the Boltzmann constant.
The logarithm in (\ref{Boltzmann_S}) is due to the fact that the entropy is
additive with respect to extensive quantities such as volume, whereas the
number of microstates is multiplicative. We will use natural logarithms
throughout, although the base of the logarithm is not important as long as
one sticks to the same choice in a calculation; a change of the logarithm
base entails just a rescaling of the entropy, i.e., a change of units.

Gibbs generalized Boltzmann's formula to systems in which the constituent
particles may fluctuate among states with different energies. This more
general scenario is made possible by setting the system in contact with a
thermal reservoir; the energy is then conserved in the composite system.
Gibbs formula for a discrete set of microstates is%
\begin{equation}
S=-k_{B}\sum \limits_{i=1}^{W}p_{i}\ln p_{i},  \label{Gibbs_S}
\end{equation}
where $p_{i}$ is the probability (i.e., the asymptotic fraction of time)
that the system is in the microstate $i$ with energy $\varepsilon_{i}$, $%
1\leq i\leq W<\infty$. Therefore 
\begin{equation}
\sum \limits_{i=1}^{W}p_{i}\varepsilon_{i}=U,  \label{canonical}
\end{equation}
where $U$ is the internal energy of the system. If $p_{i}=0$ one sets $%
0\cdot\ln0:=0$ in (\ref{Gibbs_S}).

Note that Gibbs' entropy abridges to Boltzmann's entropy for equiprobable
microstates, the so-called \textit{microcanonical ensemble}. The set of
microstates consistent with the restriction (\ref{canonical}) (along with $%
\sum_{i=1}^{W}p_{i}=1$) is called the \textit{canonical ensemble}. Regarding
the microcanonical ensemble as the actual probability distribution amounts
to Boltzmann's \textit{Ergodensatz} (ergodic hypothesis): the trajectory of
a closed system in the phase space is dense, hence the time spent in a
region is proportional to its volume. In general, the ergodic hypothesis is
not true in the physical systems.

In turn, the Tsallis entropy \cite{Tsallis1988}, 
\begin{equation}
S_{q}=\frac{k_{B}}{1-q}\left( \sum\limits_{i=1}^{W}p_{i}^{q}-1\right) ,
\label{Tsallis}
\end{equation}%
($q>0$, $q\neq 1$) generalizes the Boltzmann-Gibbs entropy $S$, Equation (%
\ref{Gibbs_S}), in the sense that $\lim_{q\rightarrow 1}S_{q}=S$. Let us
mention at this point that Havrda and Charv\'{a}t had previously introduced
in \cite{HavrdaCharvat1967} a family of entropies that differs from (\ref%
{Tsallis}) only in a factor that depends on $q$. This is the reason for the
name Havrda-Charv\'{a}t-Tsallis entropy used by some authors.

A fundamental property of (\ref{Gibbs_S}) and (\ref{Tsallis}) is concavity
with respect to the variables $p_{i}$, $1\leq i\leq W$. The concavity of
entropy and, in general, of all thermodynamical potentials (such as internal
energy, free energy, and chemical potential) guarantees the stability of
macroscopic bodies.

\subsection{Entropy in classical information theory}

In 1948 the word entropy appeared as well in the foundational papers of
Shannon on information theory, coding theory and cryptography \cite%
{Shannon1948}. According to \cite{Denbigh1990}, it was von Neumann who,
aware that Shannon's measure of information was formally the same as the
Boltzmann-Gibbs entropy, proposed him to called it entropy. Indeed, let $%
(\Omega ,\mathcal{B},\mu )$ be a probability space and $X$ a random variable
on $(\Omega ,\mathcal{B},\mu )$ with outcomes in a finite set $\Gamma
\subset \mathbb{R}$. This means that (i) $\Omega $ is a non-empty set, (ii) $%
\mathcal{B}$ is a $\sigma $-algebra of subsets of $\Omega $, (iii) $\mu $ is
a measure on the measurable space $(\Omega ,\mathcal{B})$, and (iv) $%
X:\Omega \rightarrow \Gamma $ is a measurable map, $\Gamma $ being endowed
with the discrete $\sigma $-algebra (i.e., all subsets of $\Gamma $ are
measurable). The probability mass function of $X$ is given by $p(x)=\mu
(X=x):=\mu (X^{-1}(x))$, $x\in \Gamma $.

\begin{definition}
\label{Def2.1} The (Shannon) \emph{entropy} of $X$ is 
\begin{equation}
H(X)=-\sum_{x\in\Gamma}p(x)\ln p(x).  \label{Shannon}
\end{equation}
\end{definition}

Therefore, $H(X)$ can be viewed as the expected value of the random variable 
$I(X)=-\ln p(X)$, called the information function, and interpreted as the
average information gained by knowing the outcome of the random variable $X$%
. An alternative interpretation is that $H(X)$ measures the uncertainty
about the outcome of $X$; see \cite[Sect. 1.3]{Ash1990} for other
interpretations. Note that $H(X)$ is actually a function of $p(x)$. The
range $\Gamma $ of $X$ is usually called \textit{alphabet} in information
theory. The entropy of random variables with $\Gamma =\mathbb{R}$ will be
considered in Sect. 4.2 (see Equation (\ref{H(rho)})).

More generally, let $X_{1},...,X_{n}$ be random variables on a common
probability space $(\Omega ,\mathcal{B},\mu )$ with finite alphabets $\Gamma
_{1},...,\Gamma _{n}$, respectively, and joint probability distribution 
\begin{equation*}
p(x_{1}^{n})=\mu (X_{1}=x_{1},..,X_{n}=x_{n})=\mu \left(
X_{1}^{-1}(x_{1})\cap ...\cap X_{n}^{-1}(x_{n})\right) ,
\end{equation*}%
where $x_{1}^{n}=(x_{1},...,x_{n})\in \Gamma _{1}\times ...\times \Gamma
_{n} $. The \textit{joint entropy} $H(X_{1},...,X_{n})=:H(X_{1}^{n})$ is
defined as in (\ref{Shannon}), with $p(x)$ replaced by $p(x_{1}^{n})$.

Another basic, entropy-like quantity is the \textit{conditional entropy} of
two uni- or multivariate random variables,%
\begin{equation}
H(X\left\vert Y\right) =-\sum_{x\in \Gamma _{X}}\sum_{y\in \Gamma
_{Y}}p(x,y)\ln p(x\left\vert y\right) =H(X,Y)-H(Y),  \label{Hcond}
\end{equation}%
where $p(x\left\vert y\right) =p(x,y)/p(y)$ is the conditional probability
of $X=x$ given $Y=y$, and $\Gamma _{X}$ and $\Gamma _{Y}$ are the (finite)
alphabets of $X$ and $Y$, respectively. Their \textit{mutual information} is%
\begin{equation}
I(X;Y)=-\sum_{x\in \Gamma _{X}}\sum_{y\in \Gamma _{Y}}p(x,y)\ln \frac{p(x,y)%
}{p(x)p(y)}=H(X)-H(X\left\vert Y\right) .  \label{I(X;Y)}
\end{equation}%
From (\ref{Hcond}) and (\ref{I(X;Y)}) it follows%
\begin{equation}
I(X;Y)=H(X)+H(Y)-H(X,Y)=I(Y;X).  \label{I(X;Y)B}
\end{equation}%
Therefore, $I(X;Y)$ is the reduction in the uncertainty of $X$ due to the
knowledge of $Y$; or vice versa, the average amount of information that $x$
conveys about $y$.

\begin{lemma}
\label{Lemma2.1B} \emph{(Basic properties \cite{Cover2006}.)} Let $X$, $%
X_{i} $ $(1\leq i\leq n)$, $Y$ be uni- or multivariate random variables with
finite alphabets.

\begin{enumerate}
\item Let $\Gamma $ be the alphabet of $X$ with cardinality $\left\vert
\Gamma \right\vert $. Then, $H(X)\leq \ln \left\vert \Gamma \right\vert $,
with equality if and only if $p(x)$ is uniformly distributed, i.e., $%
p(x)=1/\left\vert \Gamma \right\vert $ for all $x\in \Gamma .$

\item $H(X_{1},...,X_{n})\leq \sum_{i=1}^{n}H(X_{i})$, with equality if and
only if the $X_{i}$ are independent.

\item $H(X\left\vert Y\right) \leq H(X).$

\item $I(X;Y)\geq 0$, with equality if and only if $X$ and $Y$ are
independent.

\item \emph{(Chain rule)} $H(X_{1},...,X_{n})=\sum_{i=1}^{n}H(X_{i}\left%
\vert X_{i-1},...,X_{1}\right) .$
\end{enumerate}
\end{lemma}

Consider an \textit{information }(or \textit{data}) \textit{source} that
outputs letters, one per unit of time, from a finite alphabet $\Gamma $.
Formally, such an information source is a discrete-time, stationary
stochastic process $\mathbf{X}=(X_{t})_{t\in \mathbb{T}}$, where $\mathbb{T}=%
\mathbb{Z}$ or $\mathbb{T}=\mathbb{N}_{0}=\{0,1,..\}$ and $X_{t}$ are
equally distributed random variables on a probability space $(\Omega ,%
\mathcal{B},\mu )$ with alphabet $\Gamma $. A realization of $\mathbf{X}$ is
then a sequence $(x_{t})_{t\in \mathbb{T}}=(X_{t}(\omega ))_{t\in \mathbb{T}%
}\in \Gamma ^{\mathbb{T}}$, $\omega \in \Omega $, called a \textit{message}.
A finite segment of a message, say, $x_{t}^{t+n-1}$ is called a \textit{word}
of length $n$, and its probability to be output is the joint probability $%
p(x_{t}^{t+n-1})$. In information theory one takes usually $\mathbb{T}=%
\mathbb{N}_{0}$ because physical data sources must be turned on at some
finite time.

\begin{definition}
\label{Def2.2} The (Shannon) \emph{entropy rate} (or just \emph{entropy})%
\textit{\ }of the data source $\mathbf{X}$ is%
\begin{equation}
h(\mathbf{X})=\lim_{n\rightarrow\infty}\frac{1}{n}H(X_{0}^{n-1})=-\lim
_{n\rightarrow\infty}\frac{1}{n}\sum_{x_{0}^{n-1}\in%
\Gamma^{n}}p(x_{0}^{n-1})\ln p(x_{0}^{n-1}).  \label{Shannon_rate}
\end{equation}
\end{definition}

From the chain rule in Lemma \ref{Lemma2.1B}(5) and the stationarity of $%
\mathbf{X}$, one obtains the useful formula \cite{Cover2006} 
\begin{equation}
h(\mathbf{X})=\lim_{n\rightarrow \infty }H(X_{i}\left\vert
X_{i-1},...,X_{1}\right) .  \label{Shannon_rate2}
\end{equation}%
Since the conditional entropies $H(X_{i}\left\vert X_{i-1},...,X_{1}\right) $
build a decreasing sequence of nonnegative numbers, the convergence of (\ref%
{Shannon_rate}) and (\ref{Shannon_rate2}) follows. $h(\mathbf{X})$ measures
the average information per time unit conveyed by the messages of the
information source $\mathbf{X}$. Typical examples of information sources are
stationary, finite-state Markov chains, with state space $\Gamma $, or a
function thereof whose domain is the state space of the chain and whose
range is a subset of $\Gamma $.

The multivariate information function $I(X_{0}^{t-1})=-\ln p(X_{0}^{t-1})$
and $h(\mathbf{X})$ are related through the so-called \textit{asymptotic
equipartition principle}, this meaning that `typical' messages have roughly
the same probability. We say that an information source $\mathbf{X}$ is 
\textit{ergodic} if the long-run relative frequency of any word in a message
converges stochastically to the probability of that word.

\begin{theorem}
\label{Thm2.3} \emph{(Shannon-MacMillan-Breiman \cite{Cover2006})} Let $%
\mathbf{X}$ be an ergodic information source. Then%
\begin{equation*}
-\lim_{n\rightarrow\infty}\frac{1}{n}\ln p(X_{0}^{n-1})=h(\mathbf{X})
\end{equation*}
in probability.
\end{theorem}

Historically, Boltzmann used the ergodic hypothesis to derive the
equipartition of energy in the kinetic theory of gases. If the random
variables $X_{n}$ are independent, then Theorem \ref{Thm2.3} is a
straightforward consequence of the law of large numbers. As often happens in
other situations, ergodicity can be viewed as a property generalizing that
law.

Besides pervading all the conceptual tools of information theory, entropy is
instrumental in most applications, especially in such important ones as the
minimal code length, maximal channel capacity, or cipher security. For
brevity we will discuss next only the compression of information, which is
at the heart of modern communications and also related to Theorem \ref%
{Thm2.3}.

Compression is any procedure that reduces the data requirements of a message
without, in principle, loosing information. Suppose that code words $%
w_{1},...,w_{M}$ of lengths $L_{1},...,L_{M}$, respectively, are assigned to
the outcomes $x_{1},...,x_{M}$ of a random variable $X$ with probabilities $%
p(x_{1}),...,p(x_{M})$. The code words are combinations of letters $%
a_{1},...,a_{D}$, usually $0,1$ ($D=2$) in modern communications (or
dot/dash in the Morse code!). Then the Huffman coding is a uniquely
decipherable code that minimizes the average code-word length $\bar{L}%
=\sum_{i=1}^{M}p(x_{i})L_{i}$, which according to the \textit{noiseless
coding theorem} is known to satisfy \cite{Ash1990} 
\begin{equation}
H(X)\leq\bar{L}<H(X)+1\text{,}  \label{C1_CodingThm}
\end{equation}
the logarithm of $H(X)$ being taken to base $D$. This result allows to
interpret the entropy $H(X)$ as a lower bound on the average compression of
the symbols $x_{i}$ with relative frequencies $p(x_{i})$, $1\leq i\leq M$.

There are also algorithms for lossless data compression which, unlike the
Huffman coding, achieve the entropy bound (\ref{C1_CodingThm})
asymptotically with no prior knowledge of the probabilities $p(x_{i})$.
Examples of such `universal compressors' are different algorithms due to A.
Lempel and J. Ziv based on pattern matching \cite{Cover2006,Amigo2010A}.

Finally, let us mention that the \textit{R\'{e}nyi entropy} \cite{Renyi1961},%
\begin{equation*}
H_{q}(X)=\frac{1}{1-q}\ln \left( \sum_{x\in \Gamma }p(x)^{q}\right)
\end{equation*}%
where $q>0$, $q\neq 1$, generalizes the Shannon entropy in the sense that 
\begin{equation*}
H_{1}(X):=\lim_{q\rightarrow 1}H_{q}(X)=H(X).
\end{equation*}%
As a function of the parameter $q$, $H_{q}(X)$ is non-increasing. $%
H_{0}(X):=\lim_{q\rightarrow 0}H_{q}(X)=\ln \left\vert \Gamma \right\vert $
is also called the \textit{Hartley entropy}, $H_{2}(X)=-\ln \sum_{x\in
\Gamma }p(x)^{2}$ the \textit{collision or quadratic entropy}, and $%
H_{\infty }:=\lim_{q\rightarrow \infty }H_{q}(X)=-\ln (\max_{x\in \Gamma
}p(x))$ the \textit{min-entropy}. The R\'{e}nyi entropy, especially $%
H_{2}(X) $, has been successfully applied in information theory (see \cite%
{Principe2010} and references therein).

\subsection{Entropy in measure-preserving dynamical systems}

Entropy was introduced by Kolmogorov in ergodic theory in 1958 \cite%
{Kolmogorov1958} as a metric invariant for Bernoulli shifts. Kolmogorov's
proposal was inspired by Shannon's entropy, whose work on information theory
was well known to him. This way he was able to show that the $(\frac{1}{2},%
\frac{1}{2})$-Bernoulli shift and the $(\frac{1}{3},\frac{1}{3},\frac {1}{3}%
) $-Bernoulli shift were not metrically isomorphic because the entropy of
the former is $\ln2$ while the entropy of the latter is $\ln3$ (see below
for details). He asked then if entropy is a complete isomorphism invariant
for the Bernoulli shifts, i.e., if two Bernoulli shifts with the same
entropy are isomorphic. One year later, Sinai generalized the Kolmogorov
invariant to general measure-preserving dynamical systems \cite{Sinai1959A}
and this is the definition of entropy used ever since (also called \textit{%
metric entropy}, \textit{measure-theoretic entropy}, or \textit{%
Kolmogorov-Sinai entropy}). The question on the completeness of entropy for
Bernoulli shifts was answered affirmatively by Ornstein in 1969 \cite%
{Ornstein1970A,Ornstein1970B}. Let us point out that historically (before
the advent of chaos theory) only invertible transformations (i.e.,
automorphisms) were considered. In fact, the metric entropy is not a
complete invariant for non-invertible Bernoulli shifts \cite{Walters2000}.

To avoid technical subtleties, we assume hereafter (as usually done in
ergodic theory) that the probability spaces considered are isomorphic to the
disjoint union of countably many point masses and an interval $[a,b]\subset 
\mathbb{R}$ endowed with the Lebesgue measure. These are called \textit{%
Lebesgue spaces} \cite{Walters2000}.

\begin{definition}
\label{Def2.4} Let $(\Omega,\mathcal{B},\mu)$ be a Lebesgue probability
space and $T:\Omega\rightarrow\Omega$ a measure-preserving transformation,
that is, $T$ is measurable and $\mu(T^{-1}B)=\mu(B)$ for all $B\in\mathcal{B}
$. We say then that $(\Omega,\mathcal{B},\mu,T)$ is a \emph{%
measure-preserving} (dynamical) \emph{system}.
\end{definition}

One says equivalently that $T$ is a $\mu$-preserving map or that $\mu$ is a $%
T$-invariant measure. The transformation $T$ introduces the `dynamics' in
the state space $\Omega$ via its iterates $T^{n}$, $n$ being interpreted as
discrete time. The sequence $(T^{n}(\omega))_{n\in\mathbb{T}}$ is the 
\textit{orbit} of $\omega\in\Omega$.

The prototypes of measure-preserving dynamical systems are the \textit{shift
spaces} \textit{over} $k$ \textit{symbols} $(\{1,...,k\}^{\mathbb{T}},%
\mathcal{C},m,\sigma)$. Here $\mathcal{C}$ is the product $\sigma$-algebra
generated by the cylinder sets (i.e., sets of sequences with a finite number
of components fixed), $\sigma$ is the (left-)\textit{shift transformation} $%
(x_{t})_{t\in\mathbb{T}}\mapsto$ $(x_{t+1})_{t\in\mathbb{T}}$, and $m$ is a
probability measure that makes $\sigma$ measure-preserving. If $\mathbb{T}=%
\mathbb{Z}$ the shift is called two-sided and $\sigma$ is invertible, $%
\sigma^{-1}$ being the \textit{right shift} $(x_{t})_{t\in\mathbb{Z}}\mapsto$
$(x_{t-1})_{t\in\mathbb{Z}}$; if $\mathbb{T}=\mathbb{N}_{0}$ the shift is
called one-sided and $\sigma$ is $k$-to-$1$. Let (i) $\mathbf{p}%
=(p_{1},...,p_{k})$ be a probability vector, (ii) $P=(p_{ij})_{1\leq i,j\leq
k}$ a stochastic matrix (i.e., $p_{ij}\geq0$, $\sum%
\nolimits_{j=1}^{k}p_{ij}=1$) such that $\sum%
\nolimits_{i=1}^{k}p_{i}p_{ij}=p_{j}$, and (iii) $[i]_{t_{0}}:=%
\{(x)_{t}:x_{t_{0}}=i\}$, $1\leq i\leq k$, the $t_{0}$-time cylinder sets.
Then the $\sigma$-invariant measure $m$ of the $(\mathbf{p},P)$-\textit{%
Markov shift} is defined by $m([i]_{t_{0}})=p_{i}$, for all $t_{0}\in\mathbb{%
T}$, and extended to other cylinder sets via $m([i_{0}]_{t_{0}}\cap\lbrack
i_{1}]_{t_{1}}\cap...\cap\lbrack
i_{j}]_{t_{j}})=p_{i_{0}}p_{i_{0}i_{1}}...p_{i_{j-1}i_{j}}$. $P$ is called
the \textit{transition probability matrix} of the Markov shift. The $\mathbf{%
p}$-\textit{Bernoulli shift} is defined by taking $p_{ij}=p_{j}$, hence $%
m([i_{0}]_{t_{0}}\cap\lbrack i_{1}]_{t_{1}}\cap...\cap\lbrack
i_{j}]_{t_{j}})=p_{i_{0}}p_{i_{1}}...p_{i_{j}}$. For more general $\sigma$%
-invariant measures, see \cite[Sect. 1.1]{Walters2000}.

Given a finite partition $\alpha =\{A_{1},...,A_{k}\}$ of $\Omega $, define
the entropy of $\alpha $ as 
\begin{equation}
H_{\mu }(\alpha )=-\sum_{i=1}^{k}\mu (A_{i})\ln \mu (A_{i}).
\label{H(alpha)}
\end{equation}%
If $\beta =\{B_{1},...,B_{l}\}$ is another partition of $\Omega $, let $%
\alpha \vee \beta $ denote the \textit{least common refinement} of $\alpha $
and $\beta $, i.e., $\alpha \vee \beta =\{A_{i}\cap B_{j}:1\leq i\leq
k,1\leq j\leq l\}$. Least common refinements of more than two partitions are
defined recursively. Sinai's definition of the entropy of $(\Omega ,\mathcal{%
B},\mu ,T)$ proceeds by `coarse-graining' the state space $\Omega $ with a
partition $\alpha $ and refining $\alpha $ by the repeated action of $T$.
Specifically, define the maps $X_{t}^{\alpha }:\Omega \rightarrow
\{1,...,k\} $ by%
\begin{equation}
X_{t}^{\alpha }(\omega )=i\;\text{if\ }T^{t}(\omega )\in A_{i},\;t\in 
\mathbb{T}.\text{\ }  \label{symboldynam}
\end{equation}%
Then one can prove \cite[Sect. 1.1.3]{Amigo2010A} that $\mathbf{X}^{\alpha
}=(X_{t}^{\alpha })_{t\in \mathbb{T}}$ is a stationary random process on $%
(\Omega ,\mathcal{B},\mu )$ with alphabet $\{1,...,k\}$, called the \textit{%
symbolic dynamics of} $T$ \textit{with respect to the partition} $\alpha $.
The \textit{entropy of} $T$ \textit{with respect to the partition} $\alpha $
is then defined as the Shannon entropy of the information source $\mathbf{X}%
^{\alpha }$,%
\begin{equation}
h_{\mu }(T,\alpha )=h(\mathbf{X}^{\alpha }).  \label{h(T,a)0}
\end{equation}%
Since%
\begin{equation*}
\mu \left( X_{0}^{\alpha }=i_{0},...,X_{n-1}^{\alpha }=i_{n-1}\right) =\mu
\left( A_{i_{0}}\cap ...\cap T^{-(n-1)}A_{i_{n-1}}\right) ,
\end{equation*}%
it follows%
\begin{equation}
h_{\mu }(T,\alpha )=\lim_{n\rightarrow \infty }\frac{1}{n}H(\alpha ^{(n)}),
\label{h(T,alpha)}
\end{equation}%
where%
\begin{equation}
\alpha ^{(n)}=\alpha \vee T^{-1}\alpha \vee ...\vee T^{-n+1}\alpha
=:\bigvee\limits_{i=0}^{n-1}T^{-i}\alpha .  \label{alpha*(n)}
\end{equation}

\begin{definition}
\label{Def2.5} The \emph{(metric)} \emph{entropy} of $(\Omega,\mathcal{B}%
,\mu,T)$, or just the entropy of $T$ if the underlying measure-preserving
system is clear from the context, is 
\begin{equation}
h_{\mu}(T)=\sup_{\alpha}h_{\mu}(T,\alpha),  \label{h(T)}
\end{equation}
where the supremum is taken over all finite partitions $\alpha$ of $\Omega$.
\end{definition}

By definition, $0\leq h_{\mu }(T)\leq \infty $. Furthermore, $h_{\mu }(T)$
describes an intrinsic property of $(\Omega ,\mathcal{B},\mu ,T)$ that, most
importantly, is invariant under isomorphisms of measure-preserving dynamical
systems. A sketch of a proof is as follows. We say that $(\tilde{\Omega},%
\mathcal{\tilde{B}},\tilde{\mu},\tilde{T})$ is a \textit{factor} of $(\Omega
,\mathcal{B},\mu ,T)$ if there is a measurable map $\phi :\Omega \rightarrow 
\tilde{\Omega}$ such $\phi (T(\omega ))=\tilde{T}(\phi (\omega ))$ for $\mu $%
-almost all $\omega \in \Omega $. Then it is easily shown that $h_{\tilde{\mu%
}}(\tilde{T})\leq $ $h_{\mu }(T)$. Therefore, if $(\Omega ,\mathcal{B},\mu
,T)$ and $(\tilde{\Omega},\mathcal{\tilde{B}},\tilde{\mu},\tilde{T})$ are
(metrically) \textit{isomorphic} (i.e., each system is a factor of the
other), the invariance of the entropy follows.

A system $(\Omega,\mathcal{B},\mu,T)$ is called \textit{ergodic} if $%
T^{-1}B=B$ (or rather $\mu(T^{-1}B\bigtriangleup B)=0$), $B\in\mathcal{B}$,
implies $\mu(B)=0$ or $\mu(B)=1$. One also says that $T$ is ergodic with
respect to $\mu$, or that $\mu$ is an ergodic measure. Although this modern
definition is a far cry from Boltzmann's ergodicity hypothesis (a property
called now minimality), it is remarkable that the equality of time averages
and phase space averages of observables follows as well from the celebrated
Birkhoff's ergodic theorem \cite{Birkhoff1931}. As a useful example, the $(%
\mathbf{p},P)$-Markov shift (either one-sided or two-sided) is ergodic if
and only if $P$ is \textit{irreducible} (i.e., for all $i,j$, there exists $%
n>0$ with $p_{ij}^{n}>0$, where $p_{ij}^{n}$ is the $(i,j)$-entry of the
matrix $P^{n}$). In particular, Bernoulli shifts are ergodic.

A finite partition $\gamma=(G_{1},...,G_{k})$ of $\Omega$ is called a 
\textit{generating partition}, or just a \textit{generator} of $(\Omega ,%
\mathcal{B},\mu,T)$ if $\mathcal{B}$ is the smallest $\sigma$-algebra
containing all $T^{-n}G_{i_{n}}$, $n\in\mathbb{T}$, $1\leq i_{n}\leq k$.

\begin{theorem}
\label{Thm2.7} \emph{(The Kolmogorov-Sinai theorem \cite{Sinai1959A})} If $%
\gamma$ is a generator of the dynamical system $(\Omega,\mathcal{B},\mu,T)$,
then 
\begin{equation}
h_{\mu}(T)=h_{\mu}(T,\gamma).  \label{KS_Thm}
\end{equation}
\end{theorem}

This theorem was first proved by Kolmogorov for Bernoulli shifts using the
fact that the zero-time cylinders $\{[1]_{0},...,[k]_{0}\}$ build a
generating partition for (both the one-sided and the two-sided) $%
(p_{1},...,p_{k})$-Bernoulli shift $\sigma$. Hence,%
\begin{equation}
h_{m}(\sigma)=-\sum_{i=1}^{k}p_{i}\ln p_{i}.  \label{h_Bernoulli}
\end{equation}
In particular, the entropy of the $(\frac{1}{k},...,\frac{1}{k})$-Bernoulli
shift is $\ln k$.

Upon setting $\alpha =\gamma $ in (\ref{h(T,a)0}) one obtains 
\begin{equation}
h_{\mu }(T)=h(\mathbf{X}^{\gamma })  \label{KS=Xgamma}
\end{equation}%
in virtue of Theorem \ref{Thm2.7}. Thus generators allow to establish a
bridge between the Kolmogorov-Sinai and Shannon entropies via symbolic
dynamics. From a practical viewpoint, generators provide a method to compute
entropy via Equation (\ref{KS_Thm}). There are a number of theorems
guaranteeing the existence of generators. Thus, Krieger's theorem states
that if $T$ is ergodic and invertible with $h_{\mu }(T)<\infty $, then $T$
has a generator \cite{Krieger1970}. Although the proof is non-constructive,
Smorodinsky \cite{Smorodinsky1971} and Denker \cite{Denker1974} provided
methods to construct a generator for ergodic and aperiodic invertible maps.
Denker's construction could even be extended by Grillenberger \cite%
{Denker1976} to the nonergodic case. The existence of generators for
non-invertible maps was proved by Kowalski under different assumptions \cite%
{Kowalski1984,Kowalski1995}. However, as far as we know, the construction of
such generators remains an open problem as of this writing.

In any case, one can estimate $h_{\mu }(T)$ by taking ever finer partitions
in $h_{\mu }(T,\alpha )$. Indeed, write $(\alpha _{n})_{n\in \mathbb{N}%
}\nearrow \alpha $ if partition $\alpha _{n+1}$ is a refinement of $\alpha
_{n}$ for all $n\geq 1$, and $\alpha $ is the coarsest partition finer than
all partitions $\alpha _{n}$. Further, let $\varepsilon $ denote the
partition of $\Omega $ into separate points, i.e., $\varepsilon =\{\{\omega
\}:\omega \in \Omega \}$.

\begin{theorem}
\label{Thm2.8} \cite[Sect. 4.6]{Walters2000} Let $(\alpha _{n})_{n\in 
\mathbb{N}}$ be a sequence of partitions of $(\Omega ,\mathcal{B},\mu )$
such that $(\alpha _{n})_{n\in \mathbb{N}}\nearrow \varepsilon $. If $%
T:\Omega \rightarrow \Omega $ is a $\mu $-preserving transformation, then $%
h_{\mu }(T)=\lim_{n\rightarrow \infty }h(T,\alpha _{n})$.
\end{theorem}

Let us make the connection between dynamical systems and symbolic dynamics
more precise by introducing the concept of sequence space model. There
exists a canonical way of attaching a shift system to any \textit{stationary}
random process. For our purposes, only discrete-time, finite state processes 
$\mathbf{X}=(X_{t})_{t\in \mathbb{T}}$ are needed. Thus, let $X_{t}$ be
random variables on a probability space $(\Omega ,\mathcal{B},\mu )$ with
alphabet (without restriction) $\{1,...,k\}$, and joint probability
distributions $p(x_{t_{1}},...,x_{t_{n}})$, where $n\geq 1$, and $%
x_{t_{1}},...,x_{t_{n}}\in \{1,...,k\}$. Stationarity means $%
p(x_{t_{1}},...,x_{t_{n}})=p(x_{t_{1}+\tau },...,x_{t_{n}+\tau })$ for every 
$n\geq 1$, $t_{1},...,t_{n}\in \mathbb{T}$, and $\tau \geq 1$. Define the
map $\Phi :\Omega \rightarrow \{1,...,k\}^{\mathbb{T}}$ by $(\Phi (\omega
))_{t}=X_{t}(\omega )$ and, as before, let $\mathcal{C}$ be the product $%
\sigma $-algebra generated by the cylinder sets, i.e., sets of points in $%
\{1,...,k\}^{\mathbb{T}}$ (one- or two-sided sequences) with a finite number
of components fixed. Furthermore, endow the measurable space $(\{1,...,k\}^{%
\mathbb{T}},\mathcal{C})$ with the transported measure $m=\mu \circ \Phi
^{-1}$, i.e.,%
\begin{equation*}
m(C)=\mu (\Phi ^{-1}C),\;C\in \mathcal{C}.
\end{equation*}%
Then the shift transformation $\sigma $ on the probability space $%
(\{1,...,k\}^{\mathbb{T}},\mathcal{C},m)$ is $m$-invariant because of the
stationarity of $\mathbf{X}$.

\begin{definition}
\label{DefSupportX}\label{Def2.8} The shift space $(\{1,...,k\}^{\mathbb{T}},%
\mathcal{C},m,\sigma )$ is called the \textit{sequence space model} or 
\textit{support} of $\mathbf{X}$.
\end{definition}

The fact that historically entropy entered into the theory of dynamical
systems through the shift systems is not casual since, as we have just seen,
these are dynamical models of stationary random processes. Thus, the $%
\mathbf{p}$-Bernoulli shift models in the sense just explained an i.i.d.
random process with probability distribution $\mathbf{p}$. More generally,
the $(\mathbf{p},P)$-Markov shift $(\{1,...,k\}^{\mathbb{T}},\mathcal{C}%
,m,\sigma )$ is the support of a $k$-state, stationary Markovian process $%
(X_{t})_{t\in \mathbb{T}}$ on a probability space $(\Omega ,\mathcal{B},\mu
) $ with stationary probability distribution $\mathbf{p}$ and transition
probability matrix $P$, i.e., $\mu
(X_{t_{1}}=i_{1},...,X_{t_{n}}=i_{n})=p_{i_{1}}p_{i_{1}i_{2}}...p_{i_{n-1}i_{n}} 
$. As with the Bernoulli shifts, it can be shown that the partition $\gamma
=\{[1]_{0},...,[k]_{0}\}$ is a generator of both one- and two-sided Markov
shifts, so 
\begin{equation}
h_{m}(\sigma )=h_{m}(\sigma ,\gamma )=-\sum_{i,j=1}^{k}p_{i}p_{ij}\ln p_{ij}.
\label{h_Markov}
\end{equation}%
Of course, Equation (\ref{h_Bernoulli}) is just a special case of Equation (%
\ref{h_Markov}).

Consider now a measure-preserving system $(\Omega,\mathcal{B},\mu,T)$ and $%
\mathbf{X}^{\gamma}$, the symbolic dynamics of $T$ with respect to a
generator $\gamma=(G_{1},...,G_{k})$. Which is the relation between the
support $(\{1,...,k\}^{\mathbb{T}},\mathcal{C},m,\sigma)$ of $\mathbf{X}%
^{\gamma}$ and $(\Omega,\mathcal{B},\mu,T)$?

\begin{theorem}
\label{Thm2.9} \emph{\cite[Sect. 4.6]{Walters2000}} If $\gamma$ is a
generator of $(\Omega,\mathcal{B},\mu,T)$, then the support of the symbolic
dynamic $\mathbf{X}^{\gamma}$ is isomorphic to $(\Omega,\mathcal{B},\mu,T)$.
\end{theorem}

Therefore (see (\ref{KS=Xgamma})),%
\begin{equation*}
h_{m}(\sigma )=h_{\mu }(T)=h(\mathbf{X}^{\gamma }).
\end{equation*}%
The isomorphism is the `\textit{coding map}' $\Phi ^{\gamma }:\Omega
\rightarrow \{1,...,k\}^{\mathbb{T}}$ defined as $(\Phi ^{\gamma }(\omega
))_{n}=i$ if $T^{n}(\omega )\in G_{i}$, i.e., $\Phi ^{\gamma }(\omega )=%
\mathbf{X}^{\gamma }(\omega )$ (see Equation (\ref{symboldynam})). This
means that there is a 1-to-1 relation between $\omega \in \Omega $ and the
realizations $\mathbf{X}^{\gamma }(\omega )$ for $\mu $-almost all $\omega $%
. We conclude that the orbits of deterministic systems with known generators
can be used to produce truly random sequences. By way of illustration, the
orbits of the \textit{tent map} $\Lambda :x\mapsto \min \{x,1-x\}$, $0\leq
x<1$, are codified to tosses of a fair coin with the generator $\gamma =\{[0,%
\frac{1}{2}),[\frac{1}{2},1)\}$, since $\Phi ^{\gamma }$ is indeed an
isomorphism between $([0,1],\mathcal{B},\lambda ,\Lambda )$ and the $(\frac{1%
}{2},\frac{1}{2})$-Bernoulli shift, where $\mathcal{B}$ and $\lambda $ stand
here for the Borel $\sigma $-algebra and the Lebesgue measure of $[0,1)$,
respectively. If $\Lambda $ is replaced by the \textit{shift map} $x\mapsto
2x$ $\mbox{mod 1}$, $0\leq x<1$, then $\Phi ^{\gamma }(x)$ is the binary
expansion of $x$. The sequence $(a_{t})_{t\in \mathbb{T}}:=(X_{t}^{\gamma
}(\omega ))_{t\in \mathbb{T}}\in \{1,...,k\}^{\mathbb{T}}$ is called the $%
\gamma $-\textit{name} of $\omega $. Then, $\omega =(\Phi ^{\gamma
})^{-1}((a_{t})_{t\in \mathbb{T}})=\cap _{t\in \mathbb{T}}T^{-t}G_{a_{t}}$.

To conclude this section, let us mention that if $\Omega $ is a compact $n$%
-dimensional manifold and $T$ a $C^{1+\delta }$ diffeomorphism preserving an
absolutely continuous measure $\mu $, then $h_{\mu }(f)$ is related to the
(strictly positive) Lyapunov exponents $\chi _{i}$ (counting multiplicities)
via the celebrated Pesin's formula%
\begin{equation}
h_{\lambda }(T)=\int_{\Omega }\sum_{i=1}^{n}(\chi _{i})_{+}\mathrm{d}\mu ,
\label{Pesin}
\end{equation}%
where $(\cdot )_{+}$ stands for positive part \cite{Pesin1977}. If $T$ is
ergodic, then $\chi _{i}$ are constant $\mu $-almost everywhere and the
integration in (\ref{Pesin}) can be dropped. A special case of this theorem
is of historical relevance, namely, Sinai's formula for the entropy of
automorphisms of the $n$-dimensional torus \cite{Sinai1959B}. Let $\Omega =%
\mathbb{R}^{n}/\mathbb{Z}^{n}$ and $T_{M}:\omega \mapsto M\omega $ $(%
\mbox{mod 1})$, where $M$ is an $n\times n$ matrix with integer entries and
determinant of absolute value one. Then the entropy of $T$ with respect to
the Lebesgue measure $\lambda $ is given by the eigenvalues $\Lambda
_{1},...,\Lambda _{n}$ of $M$ as%
\begin{equation}
h_{\lambda }(T_{M})=\sum_{i=1}^{n}(\ln \left\vert \Lambda _{i}\right\vert
)_{+}.  \label{h_auto}
\end{equation}%
This result showed for the first time that not only random processes but
also `deterministic processes' can have a positive entropy. Pesin's formula
was generalized to arbitrary ergodic Borel measures by Ledrappier and Young 
\cite{Ledrappier1985}.

\subsection{Entropy in topological dynamical systems}

Let $\Omega$ be a compact metric space and $T:\Omega\rightarrow\Omega$ a
continuous transformation. We say then that $(\Omega,T)$ is a \textit{%
topological dynamical system}. The notion of entropy in $(\Omega,T)$, or
topological entropy, was introduced by Adler, Konheim, and McAndrew \cite%
{Adler1965}. Their proposal follows the definition of metric entropy but
using open covers instead of partitions. Specifically, if $\alpha
=\{A_{1},...,A_{k}\}$ is now an open cover of $\Omega$, then the \textit{%
topological entropy of} $T$ \textit{relative to} $\alpha$ is given by%
\begin{equation}
h_{top}(T,\alpha)=\lim_{n\rightarrow\infty}\frac{1}{n}\ln N(\alpha^{(n)}),
\label{h_top(T,alpha)}
\end{equation}
where, similarly to (\ref{alpha*(n)}), $\alpha^{(n)}$ is the least common
refinement of the open covers $\alpha,T^{-1}\alpha...,T^{-n+1}\alpha$ (i.e., 
$\alpha^{(n)}=\{A_{i_{0}}\cap...\cap T^{-n+1}A_{i_{n-1}}:1\leq
i_{0},...,i_{n-1}\leq k\}$), and $N(\alpha^{(n)})$ denotes the number of
sets in a finite subcover of $\alpha^{(n)}$ with smallest cardinality.

\begin{definition}
\label{Def2.10} The \textit{topological entropy} of $T$ is defined by 
\begin{equation}
h_{top}(T)=\sup_{\alpha}h_{top}(T,\alpha),  \label{h_top(T)}
\end{equation}
where $\alpha$ ranges over all open covers of $\Omega$.
\end{definition}

Later, Dinaburg \cite{Dinaburg1070} and Bowen \cite{Bowen1971} introduced a
different definition for continuous maps on compact metric spaces $(\Omega
,\rho )$ which is equivalent to Definition \ref{Def2.10} (thus, independent
of the metric $\rho $) and, moreover, can be extended to noncompact metric
spaces. In this approach, one introduces a new metric%
\begin{equation}
\rho _{k}(\omega ,\omega ^{\prime })=\max_{0\leq i\leq k-1}\rho
(T^{i}(\omega ),T^{i}(\omega ^{\prime }))  \label{rhok}
\end{equation}%
in $\Omega $ that takes into account the separation of points along their
initial orbit segments of length $k\geq 1$. Denote by $r_{k}(\epsilon )$ the
minimal number of $\epsilon $-balls with respect to $\rho _{k}$ that cover
the whole space. Then 
\begin{equation}
h_{top}(T)=\lim_{\epsilon \rightarrow 0}\underset{k\rightarrow \infty }{\lim
\sup }\frac{1}{k}\log r_{k}(\epsilon )=\lim_{\epsilon \rightarrow 0}\underset%
{k\rightarrow \infty }{\lim \inf }\frac{1}{k}\log r_{k}(\epsilon ).
\label{Bowen}
\end{equation}%
Thus, the topological entropy measures the exponential growth rate of the
number of distinguishable orbits with finite precision. For other
approaches, also when $\Omega $ is noncompact, see \cite[Sect. 7.2]%
{Walters2000}. For a survey of topological entropy in different settings of
topological dynamics, see \cite{Llibre2015}.

Topological entropy is an invariant of the equivalence between topological
dynamical systems, a notion usually called \textit{conjugacy}. Two
topological dynamical systems $(\Omega,T)$ and $(\tilde{\Omega},\tilde{T})$
are said to be (topologically) conjugate if there exists a homeomorphism $%
\phi:\Omega \rightarrow\tilde{\Omega}$ such that $\phi\circ T=\tilde{T}%
\circ\phi$, in which case $h_{top}(T)=h_{top}(\tilde{T})$.

A prototype of a \textit{topological dynamical system} $(\Omega,T)$ is the 
\textit{full shift over }$k$ \textit{symbols} $(\{1,...,k\}^{\mathbb{T}%
},\sigma)$, where $\{1,...,k\}$ is equipped with the discrete topology and $%
\{1,...,k\}^{\mathbb{T}}$ with the product topology (generated by the
cylinder sets). Indeed, $\{1,...,k\}^{\mathbb{T}}$ is then a compact,
metrisable space, and the shift $\sigma$ is a continuous transformation
(bi-continuous, hence a homeomorphism, if $\mathbb{T}=\mathbb{Z}$). More
interesting though are the \textit{subshifts}, which are obtained from a
full shift by excluding a finite or infinite set of `forbidden' words from
the state space. Therefore, subshifts are shift spaces with constrained
sequences. Let $S$ be the subset of $\{1,...,k\}^{\mathbb{T}}$ built by all
allowed words, also called a \textit{language}. The topological entropy of
the subshift $\sigma_{S}=\left. \sigma\right\vert _{S}$ is given by%
\begin{equation}
h_{top}(\sigma_{S})=\lim_{n\rightarrow\infty}\frac{1}{n}\ln\left\vert 
\mathcal{A}_{n}\right\vert ,  \label{htop(sigmaW)}
\end{equation}
where $\mathcal{A}_{n}\subset S$ is the set of \textit{allowed} words of
length $n$. If the list of forbidden words is finite, one speaks of a 
\textit{subshift of finite type} (SFT).

Two SFTs $\sigma_{S}$ and $\sigma_{S^{\prime}}$ are said to be \textit{%
almost conjugate} if there is a third SFT $\sigma_{R}$ and factor maps $\phi
:R\rightarrow S$, $\phi^{\prime}:R\rightarrow S^{\prime}$ that are 1-to-1 on
an open dense set. Topological entropy is a complete invariant of almost
topological conjugacy for aperiodic and irreducible SFTs \cite{Adler1979}.
In general, a shift space $(S,\sigma)$ is \textit{irreducible} if there is a
point $\omega\in S$ whose forward orbit $(\sigma^{n}(\omega))_{n\in \mathbb{N%
}_{0}}$ is dense in $\Omega$.

A (\textit{topological})\textit{\ Markov chain} (or 1-step SFT) is a SFT $%
(S,\sigma _{S})$ such that, given an allowed word $x_{t}...x_{t+n}$ and a
letter $x$, the concatenation $x_{t}...x_{t+n}x$ is an allowed word if and
only if $x_{t+n}x$ is an allowed word. This being the case, a Markov chain
can be ascribed a so-called transition matrix $A=(a_{ij})_{1\leq i,j\leq k}$%
, $a_{ij}\in \{0,1\}$, in such a way that $S=\{(x_{t})_{t\in \mathbb{T}%
}:a_{x_{t}x_{t+1}}=1$ for all $t\in \mathbb{T\}}$. A Markov chain is
irreducible if and only if the matrix $A$ is irreducible. Alternatively, a
Markov chain can be described by a graph $\mathcal{G}$ and its adjacency
matrix $A_{\mathcal{G}}=A$ (up to reordering of the vertices). In either
case, $h(\sigma _{A})=\ln \rho _{A}$, where $\rho _{A}$ is the spectral
radius of $A$ (i.e., the largest absolute value of the eigenvalues). In
particular, $h_{top}(\sigma )=\ln k$ for the one- and two-sided full shifts.
Note for further reference that $\ln k$ is also the metric entropy of the $(%
\frac{1}{k},...,\frac{1}{k})$-Bernoulli shift. A factor of a Markov chain is
called a \textit{sofic shift} \cite{Weiss1973}.

Subshifts were considered by Shannon in the context of constrained data
sources \cite{Shannon1948}. Indeed, due to technological feasibility or
convenience, it is sometimes necessary to encode messages to sequences that
satisfy certain constrains. For example, to ensure proper synchronization in
magnetic or optical recording it might be necessary to limit the length of
runs of 0's between two 1's when reading and recording bits \cite{Cover2006}%
. Specifically, if $\mathcal{W}_{n}$ is the set of \textit{allowed} words of
length $n$ that can be transmitted over a noiseless channel, then%
\begin{equation}
C=\lim_{n\rightarrow \infty }\frac{1}{n}\ln \left\vert \mathcal{W}%
_{n}\right\vert  \label{Capacity}
\end{equation}%
is the capacity of the `constrained channel'. By definition of noiseless
channel, $C$ is an upper bound for the entropy of any data source
transmitting over the channel. Comparison of (\ref{Capacity}) with (\ref%
{htop(sigmaW)}) shows that $C$ is the topological entropy of the subshift $%
\sigma _{S}$ with language $S=\cup _{n\geq 1}\mathcal{W}_{n}$. Shannon
proved \cite{Shannon1948} that if the data source is an irreducible
Markovian random process (of order 1), then there is an assignment of
transition probabilities which maximizes the entropy. In other words, if $%
(S,\sigma _{S})$ is an irreducible Markov chain, there exists a $\sigma _{S}$%
-invariant measure $m^{\ast }$ on $(S,\mathcal{C})$ such that the metric
entropy of $(S,\mathcal{C},m^{\ast },\sigma _{S})$ equals the topological
entropy of $(S,\sigma _{S})$. $m^{\ast }$ is called the \textit{Parry measure%
} (he proved its uniqueness in \cite{Parry1964}); $(S,\mathcal{C},m^{\ast
},\sigma _{S})$ is the support of the capacity achieving Markovian source.
The study of constrained languages evolved over time to symbolic dynamics,
nowadays a stand-alone branch of mathematics with important applications to
areas other than dynamical systems and information theory, such as formal
languages, computer science, and graph theory \cite{Lind1995}.

According to the Krylov-Bogolyubov theorem, any topological dynamical system
has at least one Borel probability $T$-invariant measure \cite%
{Hasselblatt2002}. The proof resorts to standard results of functional
analysis. Set $\mu_{n}(\omega):=\frac{1}{n}\sum_{i=0}^{n-1}\delta
_{T^{i}(\omega)}$, $\omega\in\Omega$, where $\delta_{T^{i}(\omega)}$ is the
unit point mass at $T^{i}(\omega)$. Then every $\mu_{n}(\omega)$ is a Borel
probability measure on $\Omega$, and the weak$^{\ast}$-limit of $(\mu
_{n}(\omega))_{n\in\mathbb{N}_{0}}$ is a $T$-invariant Borel probability
measure on $\Omega$, supported on the closure of the orbit $%
(T^{i}(\omega))_{i\in\mathbb{N}_{0}}$. Note that distinct periodic orbits
(if any) contribute distinct invariant measures. The existence of invariant
measures in any topological dynamical system raises the question about the
relationship between topological and metric entropy. The answer is given by
the following \textit{variational principle}.

\begin{theorem}
\label{Thm2.11} \emph{(\cite{Goodman1971})}. Let $\Omega$ be a compact
metric space and $T:\Omega\rightarrow\Omega$ a continuous map. Then 
\begin{equation*}
h_{top}(T)=\sup_{\mu}h_{\mu}(T),
\end{equation*}
where $\mu$ ranges over all Borel probability $T$-invariant measures.
\end{theorem}

The above examples of the $(\frac{1}{k},...,\frac{1}{k})$-Bernoulli measure
and the Parry measure show that the bound $h_{\mu}(T)\leq h_{top}(T)$ is
actually sharp. A $T$-invariant measure $\mu^{\ast}$ such that $%
h_{\mu^{\ast}}(T)=h_{top}(T)$ is called a \textit{measure of maximal entropy}%
. If, moreover, $\mu^{\ast}$ is unique (as the both examples just
mentioned), then $\mu^{\ast}$ is the natural choice to characterize the
dynamics.

Both the topological entropy and the variational principal are conveniently
generalized in the so-called \textit{thermodynamical formalism} \cite%
{Sinai1972}. One of the main characters in this approach is the topological
pressure \cite{Ruelle1973}. Let $C(\Omega ,\mathbb{R})$ be the space of
real-valued continuous functions of a compact metric space $\Omega $. The 
\textit{pressure of} $\mu $, a $T$-invariant Borel probability measure, is
the map $P_{\mu }(T,\cdot ):C(\Omega ,\mathbb{R})\rightarrow \mathbb{R}\cup
\{\infty \}$ defined by $P_{\mu }(T,\varphi )=h_{\mu }(T)+\int_{\Omega
}\varphi \mathrm{d}\mu $. The function $\varphi $ is sometimes called the 
\textit{potential}. Furthermore, we define the (topological) \textit{%
pressure of} $T$ with respect to $\varphi \in C(\Omega ,\mathbb{R})$ as%
\begin{equation}
P(T,\varphi )=\sup_{\mu }P_{\mu }(T,\varphi ).  \label{pressure}
\end{equation}%
An invariant measure $\mu ^{\ast }$ such that $P_{\mu ^{\ast }}(T,\varphi
)=P(T,\varphi )$ is called an \textit{equilibrium} \textit{state} for $%
\varphi $. Since $P(T,0)=h_{top}(T)$, a measure of maximal entropy is then
an equilibrium state for the potential $\varphi =0$. One of the most
important properties of $P(T,\varphi )$ is its characterization of invariant
measures:\ A Borel probability measure $\mu $ is $T$-invariant (where $%
h_{top}(T)<\infty $) if and only if $\int_{\Omega }\varphi \mathrm{d}\mu
\leq P(T,\varphi )$ for all $\varphi \in C(\Omega ,\mathbb{R})$.

Beside Markov chains, there are few systems for which exact formulas or good
algorithms to compute the topological entropy are known. One of the
exceptions are the \textit{multimodal maps}, which are relevant in the study
of low dimensional chaos \cite{Milnor1988B}. See also \cite%
{Amigo2014B,Amigo2015B} for recursive algorithms based on the expression%
\begin{equation*}
h(f)=\lim_{n\to\infty}\frac{1}{n}\ln\ell_{n}
\end{equation*}
where $\ell_{n}$ is the number of maximal monotonicity intervals (`lap
number') of $f^{n}$ \cite{Misiurewicz1980}.

Topological entropy was extended to non-autonomous dynamical systems in \cite%
{Kolyada1996}. In the case of switching systems, a different approach was
used in \cite{Amigo2013A,Amigo2013B}. Other extensions, this time to
set-valued mappings, can be found in \cite{Carrasco2015}.

\subsection{Entropy of group actions}

Let $(\Omega ,\mathcal{B},\mu )$ be a Lebesgue probability space and $%
T:\Omega \rightarrow \Omega $ a $\mu $-preserving transformation. The set $%
\mathcal{M}_{\mu }(\Omega )$ of $\mu $-preserving transformations $T:\Omega
\rightarrow \Omega $ is a group if $T$ is invertible or, otherwise, a
semigroup under composition.

\begin{definition}
\label{Def_MPGA}Let $G$ be a group or semigroup. A \emph{measure-preserving}%
\textit{\ }$G$\textit{\ }\emph{action on} $(\Omega,\mathcal{B},\mu)$ is a
representation of $G$ by means of $\mu$-preserving transformations $%
T:\Omega\rightarrow\Omega$.
\end{definition}

This means that there exists a homomorphism $\psi:G\rightarrow\mathcal{M}%
_{\mu}(\Omega)$, that is, $\psi(g)=T_{g}$ and $\psi(h)=T_{h}$ implies $%
\psi(gh)=T_{g}\circ T_{h}$ for every $g,h\in G$, where $(T_{g}\circ
T_{h})(\omega):=T_{g}(T_{h}(\omega))$. If $G$ is a group with identity
element $e$, then $\psi(e)$ is the identity transformation, and $%
T_{g^{-1}}=T_{g}^{-1}$. Sometimes one says that $T_{g}$ is a
(measure-preserving) action of $G$ on $\Omega$.

A $\mathbb{Z}$ or $\mathbb{N}_{0}$ action on $(\Omega,\mathcal{B},\mu)$
amounts to the invertible or non-invertible measure-preserving dynamical
system $(\Omega,\mathcal{B},\mu,T)$, respectively, where $T$ generates $G$,
i.e., $\psi(n)=T^{n}$. For another example, let $(\Omega,\mathcal{B},\mu)$
be a compact abelian group with the Haar measure $\mu$. Furthermore, let $G$
be a subgroup of $\Omega$, and the transformation $T_{g}:\Omega\rightarrow%
\Omega$ be given by $T_{g}(\omega)=g\omega$, $g\in G$. Then $g\mapsto T_{g}$
is a measure-preserving $G$ action on $(\Omega,\mathcal{B},\mu)$.

The perhaps simplest example of a measure-preserving action of an arbitrary
group $G$ is a Bernoulli action \cite{Weiss2015}. Let $K$ be a finite set
with the discrete topology, and $\Omega=K^{G}:=\{\omega:G\rightarrow K\}$
with the product topology, so $\Omega$ is compact. Define a right action $%
R_{g}\omega(h)=\omega(hg)$ and a left action $L_{g}\omega(h)=\omega(g^{-1}h)$
of $G$ on $\Omega$. To construct an invariant measure for these two actions
endow $K$ with the discrete $\sigma$-algebra $\mathcal{D}$ and a probability
mass function $p$, and endow $\Omega$ with the product measure, taking for
each coordinate the probability space $(K,\mathcal{D},p)$. For $G=\mathbb{Z}$
one obtains the standard (two-sided) Bernoulli shift on $\left\vert
K\right\vert $ symbols, with $\sigma^{n}=R_{n}=L_{-n}$, $n\in\mathbb{Z}$. A
similar construction can be made for a semigroup, where now only the right
action $R_{g}$ is defined. For $G=\mathbb{N}_{0}$ one obtains the one-sided
Bernoulli shift.

$\mathbb{R}$-actions are also called \textit{flows}. A prominent example of
an $\mathbb{R}$ action is the \textit{Bernoulli flow}, i.e., a
representation $\mathbb{R}\ni x\mapsto T_{x}\in\mathcal{M}_{\mu}(\Omega)$,
such that $T_{x}$ is isomorphic to a Bernoulli shift \cite{Ornstein2013}.

Set $\mathbb{K}=\mathbb{Z}$ or $\mathbb{N}_{0}$. The generalization of
measure-preserving $\mathbb{K}$ actions to measure-preserving $\mathbb{K}%
_{1}\times...\times\mathbb{K}_{d}$ actions on $(\Omega,\mathcal{B},\mu)$ is
as follows. Let $T_{i}:\Omega\rightarrow\Omega$, $1\leq i\leq d$, be $\mu $%
-preserving (invertible if $\mathbb{K}_{i}=\mathbb{Z}$) transformations such
that they commute with each other. For each $\mathbf{n}=(n_{1},...,n_{d})\in%
\mathbb{K}_{1}\times...\times\mathbb{K}_{d}$, define the action $T^{\mathbf{n%
}}:=T_{1}^{n_{1}}\circ...\circ T_{d}^{n_{d}}$ of $\mathbb{K}%
_{1}\times...\times\mathbb{K}_{d}$ on $\Omega$.

The definition of metric entropy for $\mathbb{K}_{1}\times ...\times \mathbb{%
K}_{d}$ actions is a straightforward generalization of the definition of
metric entropy. Indeed, let $\alpha $ be a finite partition of $\Omega $,
and set%
\begin{equation*}
\alpha ^{(n)}=\bigvee\limits_{\mathbf{n}\in \{0,1,...,n-1\}^{d}}T^{-\mathbf{n%
}}\alpha ,
\end{equation*}%
and%
\begin{equation}
h_{\mu ,d}(T,\alpha )=\lim_{n\rightarrow \infty }\frac{1}{n^{d}}H_{\mu
}(\alpha ^{(n)}).  \label{action2}
\end{equation}%
Then,%
\begin{equation}
h_{\mu ,d}(T)=\sup_{\alpha }h_{\mu ,d}(T,\alpha )  \label{action3}
\end{equation}%
is the $d$-\textit{dimensional entropy of the measure-preserving }$\mathbb{K}%
_{1}\times ...\times \mathbb{K}_{d}$ \textit{action} on $(\Omega ,\mathcal{B}%
,\mu )$. For $d=1$ we recover the definition of the Kolmogorov-Sinai entropy.

The extension of metric entropy from $\mathbb{Z}$ to more general groups has
been done in two major steps. In a first step, Ornstein and Weiss \cite%
{Ornstein1987} extended the necessary theoretical framework and fundamental
results (such as Kolmogorov's theorem for Bernoulli shifts, and Ornstein's
isomorphism theorem) to\ the \textit{amenable groups}. These are groups to
which Riesz's proof of von Neumann ergodic theorem carries over.
Incidentally, the amenable groups were introduced and studied by von Neumann
in connection with the Banach-Tarski paradox. They are precisely those
locally compact groups for which no paradoxical decomposition exists. For
the purposes of ergodic theory, a group $G$ is amenable if it has a sequence
of finite sets $F_{n}$, called a \textit{F\o lner sequence}, such that%
\begin{equation*}
\lim_{n\rightarrow \infty }\frac{\left\vert gF_{n}\Delta F_{n}\right\vert }{%
\left\vert F_{n}\right\vert }=0
\end{equation*}%
for any $g\in G$. F\o lner sequences\textit{\ }are actually the handle that
makes amenable groups amenable to the concepts and tools of ergodic theory,
in particular to entropy \cite{Rudolph2000}. It can be shown that abelian
(e.g., $\mathbb{Z}^{d}$), nilpotent and solvable groups are amenable.

In a second step, Lewis Bowen proved that an entropy theory can be also
developed for the much larger class of \textit{sofic groups}, introduced by
Gromov \cite{Bowen2010A,Bowen2010B}. The interested reader is referred to
the excellent review \cite{Weiss2015}.

Likewise, topological dynamical systems can be also generalized to \textit{%
continuous actions} of groups as follows. Recall that a topological space is
called $\sigma$-compact if it is a countable union of compact sets. Given a $%
\sigma$-compact space $\Omega$, the set of continuous transformations $%
T:\Omega\rightarrow\Omega$ is a group (if $T$ is invertible) or a semigroup
(otherwise) under composition. We focus on the group of homeomorphisms.

\begin{definition}
\cite[Definition 8.1]{Einsiedler2011} Let $G$ be $\sigma$-compact metric
group and $\Omega$ a $\sigma$-compact metric space. A \emph{continuous}%
\textit{\ }$G$\textit{\ }\emph{action on} $\Omega$ is a homomorphism from $G$
to the group of homeomorphisms of $\Omega$, $g\mapsto T_{g}$, such that the
map $G\times\Omega\rightarrow\Omega$ defined by $(g,\omega)\mapsto
T_{g}(\omega)$ is continuous. A Borel measure $\mu$ on $\Omega$ is invariant
under $G$ if $\mu\circ T_{g}^{-1}=\mu$ for all $g\in G$, where $(\mu\circ
T_{g}^{-1})(B):=\mu(T_{g}^{-1}B)$ for all Borel sets $B\subset\Omega$.
\end{definition}

The usual definition of topological entropy for $G=\mathbb{Z}$ based on
spanning and separating sets (see Sect. 2.4) carries over smoothly to
amenable groups with the help of F\o lner sequences\textit{.} The same
happens with the Krylov-Bogolyubov theorem.

\begin{theorem}
\cite[Theorem 8.10]{Einsiedler2011} If a locally compact group $G$ is
amenable, then every continuous $G$ action $T_{g}$ on a compact metric space
has an invariant probability measure.
\end{theorem}

It is therefore very satisfactory that the topological entropy of the action
of an amenable group $G$ on a compact metric space $\Omega$ and the
entropies of the same action viewed as a measure-preserving one on $\Omega$
are related by a variational principle. The same happens with sofic groups.
See \cite{Weiss2015} for further details and references.

\section{Entropy is a very natural concept}

It is not surprising that, owing to the paramount role of entropy in
information theory, Shannon pondered over its uniqueness from the very
beginning. Indeed, he proved in Appendix 2 of his foundational paper \cite%
{Shannon1948} that the only function $H$ of a probability mass distribution $%
\{p_{1},...,p_{n}\}$ satisfying just three properties (see below) is of the
form $H=-K\sum_{i=1}^{n}p_{i}\ln p_{i}$, where $K$ is a positive constant.

\subsection{Entropy from an axiomatic viewpoint}

There exist various axiomatic characterizations of entropy showing that it
is a very natural construct, in particular under its interpretation as a
measure of information or uncertainty. We list here only a few of them for
the Shannon, R\'{e}nyi and Tsallis entropies.

If $n\in $ ${\mathbb{N}}$, set 
\begin{equation*}
P_{n}=\left\{ (p_{1},p_{2},\ldots ,p_{n})\in \lbrack 0,1]^{n}\ :\
\sum\limits_{i=1}^{n}p_{i}=1\right\}
\end{equation*}%
and $\mathcal{P}=\bigcup_{n=1}^{\infty }P_{n}$. For any map $H:\mathcal{P}%
\rightarrow \lbrack 0,\infty )$, consider the following properties:

\begin{enumerate}
\item[(P1)] \textsc{Continuity:} $H$ is continuous on $P_{n}$ for each $n$.

\item[(P2)] \textsc{Symmetry:} For each $(p_{1},p_{2},\ldots ,p_{n})\in
P_{n} $ and each permutation $\pi $ of $\{0,1,\ldots ,n\}$, 
\begin{equation*}
H(p_{1},p_{2},\ldots ,p_{n})=H(p_{\pi (1)},p_{\pi (2)},\ldots ,p_{\pi (n)}).
\end{equation*}

\item[(P3)] \textsc{Monotonicity:} If $n<m$, then 
\begin{align*}
H\left( \frac{1}{n},\frac{1}{n},\ldots,\frac{1}{n}\right) <H\left( \frac {1}{%
m},\frac{1}{m},\ldots,\frac{1}{m}\right) .
\end{align*}

\item[(P4)] \textsc{Maximality:} For each $(p_{1},p_{2},\ldots,p_{n})\in
P_{n}$, 
\begin{align*}
H(p_{1},p_{2},\ldots,p_{n})\leq H\left( \frac{1}{n},\frac{1}{n},\ldots ,%
\frac{1}{n}\right) .
\end{align*}

\item[(P5)] \textsc{Expansibility:} For each $(p_{1},p_{2},\ldots ,p_{n})\in
P_{n}$ and $i\in \{1,\ldots ,n-1\}$,%
\begin{eqnarray*}
H(p_{1},p_{2},\ldots ,p_{n}) &=&H(0,p_{1},\ldots ,p_{n}) \\
&=&H(p_{1},\ldots ,p_{i},0,p_{i+1},\ldots ,p_{n})=H(p_{1},\ldots ,p_{n},0).
\end{eqnarray*}

\item[(P6)] $a$\textsc{-Additivity:} For each $(p_{1},p_{2},\ldots
,p_{m})\in P_{m}$, $(q_{1},q_{2},\ldots ,q_{n})\in \nolinebreak P_{n}$ and $%
a>0$, 
\begin{align*}
& H(p_{1}q_{1},\ldots ,p_{1}q_{n},p_{2}q_{1},\ldots ,p_{2}q_{n},\ldots
,p_{m}q_{1},\ldots ,p_{m}q_{n})\hspace*{3cm} \\
& =H(p_{1},\ldots ,p_{m})+H(q_{1},\ldots ,q_{n})+(1-a)H(p_{1},\ldots
,p_{m})H(q_{1},\ldots ,q_{n}).
\end{align*}

\item[(P7)] \textsc{Strong }$a$\textsc{-additivity:} For each $%
(p_{11},\ldots ,p_{1n},p_{21},\ldots ,p_{2n},\ldots ,p_{m1},\ldots ,
p_{mn})\in P_{mn}$, $p_{i\cdot }:=\sum_{j=1}^{n}p_{ij}$, $p_{\cdot
j}:=\sum_{i=1}^{m}p_{ij}$ and $a>0$, 
\begin{align}
& H(p_{11},\ldots ,p_{1n},p_{21},\ldots p_{2n},\ldots ,p_{m1},\ldots ,p_{mn})%
\hspace*{3cm}  \label{add} \\
& =H(p_{1\cdot },p_{2\cdot },\ldots ,p_{m\cdot })+\sum_{i=1}^{m}p_{i\cdot
}^{a}H\left( \frac{p_{i1}}{p_{i\cdot }},\frac{p_{i2}}{p_{i\cdot }},\ldots ,%
\frac{p_{in}}{p_{i\cdot }}\right) .  \notag
\end{align}

\item[(P8)] $a$\textsc{-Recursivity:} For each $(p_{1},p_{2},\ldots
,p_{n})\in P_{n},\ n>2$, and $a>0$, 
\begin{equation*}
H(p_{1},p_{2},\ldots ,p_{n})=H(p_{1}+p_{2},p_{3},\ldots
,p_{n})+(p_{1}+p_{2})^{a}\,H\left( \frac{p_{1}}{p_{1}+p_{2}},\frac{p_{2}}{%
p_{1}+p_{2}}\right) .
\end{equation*}
\end{enumerate}

Bearing in mind that $H$ shall measure the information or uncertainty
contained in a probability distribution $(p_{1},p_{2},\ldots,p_{n})$,
properties (P1)-(P8) are well interpretable.

In particular, if $p_{1},p_{2},\ldots ,p_{m}$ and $q_{1},q_{2},\ldots ,q_{n}$
are the probabilities of the outcomes $i$ and $j$ of two given
finitely-valued random variables $X$ and $Y$, respectively, then property
(P6) for $a=1$ expresses that, in case that $X$ and $Y$ are stochastically
independent, then the information contained in $X$ and $Y$ adds to the
information contained in the random vector $(X,Y)$. For $a=1$, property (P6)
is simply called \emph{additivity}.

Moreover, (P7) for $a=1$ can be interpreted as follows: If $%
p_{ij},i=1,2,\ldots ,m,j=1,2,\ldots ,n,$ are the joint probabilities of the
outcomes $i$ and $j$ of two random variables $X$ and $Y$, respectively, then
the information contained in the random vector $(X,Y)$ is the sum of the
information contained in $X$ and the mean information of $Y$ given the
outcome of $X$, that is, $H(X,Y)=H(X)+H(Y\left\vert X\right) $ in the
notation of Section 2.2 (see Equation (\ref{Hcond})).

Obviously, for $a=1$ (P7) implies (P6).

\begin{theorem}
\label{Thm3.1}Let $H:\mathcal{P}\rightarrow \lbrack 0,\infty )$ and, for $%
q>0 $ with $q\neq 1$, let the functions $f,f_{q},g_{q}:\mathcal{P}%
\rightarrow \lbrack 0,\infty )$ be defined by 
\begin{align*}
f(p_{1},p_{2},\ldots p_{n})& =-\sum_{i=1}^{n}p_{i}\ln p_{i}, \\
f_{q}(p_{1},p_{2},\ldots p_{n})& =\frac{1}{1-q}\ln \left(
\sum_{i=1}^{n}p_{i}^{q}\right) , \\
g_{q}(p_{1},p_{2},\ldots ,p_{n})& =\frac{1}{1-q}\left(
\sum_{i=1}^{n}p_{i}^{q}-1\right) ,
\end{align*}%
for all $(p_{1},p_{2},\ldots p_{n})\in \mathcal{P}$, $n\in {\mathbb{N}}$.
Then the following statements hold.

\begin{enumerate}
\item[(i)] $H=f$ satisfies (P1)-(P5), and also (P6)-(P8) with $a=1$.

\item[(ii)] $H=f_{q}$ satisfies (P1)-(P5), but not (P7) nor (P8) for any $%
a>0 $, and it is additive.

\item[(iii)] $H=g_{q}$ satisfies (P1)-(P5), and also (P7) and (P8) with $a=q$%
, but it is non-additive.

\item[(iv)] The following statements are equivalent:

\begin{enumerate}
\item[(a)] $H=c\,f$ for some $c\geq0$.

\item[(b)] $H$ satisfies (P1), (P4), (P5), and (P7) with $a=1$.

\item[(c)] $H$ satisfies (P1), (P2), and (P8) with $a=1$.
\end{enumerate}

\item[(v)] The following statements are equivalent:

\begin{enumerate}
\item[(a)] $H=c\,g_{q}$ for some $c\geq0$.

\item[(b)] $H$ satisfies (P1), (P2), and (P8) with $a=q$.
\end{enumerate}
\end{enumerate}
\end{theorem}

Of course, $f,f_{q},g_{q}$ correspond to the Shannon, R\'{e}nyi and Tsallis
entropies, respectively. Recall that R\'{e}nyi and Tsallis entropies for $%
q=1 $ are defined as the Shannon entropy since $f_{q}$ and $g_{q}$ converge
to $f $ as $q\rightarrow 1$. Furthermore, for fixed $q$ the R\'{e}nyi
entropy $f_{q}$ and the Tsallis entropy $g_{p}$ are monotonically related as
follows: 
\begin{equation*}
f_{q}=\frac{1}{1-q}\ln [1+(1-q)g_{q}]\text{\ \ or, equivalently,\ \ }g_{q}=%
\frac{1}{1-q}(e^{(1-q)f_{q}}-1).
\end{equation*}

The equivalences of (a) and (b) and (a) and (c) in Theorem \ref{Thm3.1}(iv),
both characterizing the Shannon entropy $f$, are well known results of
Khinchin \cite{Khinchin1957} and Faddeev \cite{Faddeev1956}, respectively.
Here, the properties (P7) and (P8), respectively, with $a=1$ are the
substantial ones; in their general form they play a central role for
characterizing Tsallis entropy $g_{q}$. Shannon already characterized his
entropy in \cite{Shannon1948}, by the properties (P1), (P3) and a third
property substantially being the strong $1$-additivity (additivity `for
successive choices'). Probability functionals like the R\'{e}nyi and Tsallis
entropies that satisfy the properties of continuity (P1), maximality (P4)
and expansibility (P5) are called \textit{generalized entropies}.

The characterization of the Tsallis entropy $g_{q}$ based on (P8) in Theorem %
\ref{Thm3.1}(v) was given by Furuichi \cite{Furuichi2005}. A statement
generalizing the equivalence of (a) and (b) in Theorem \ref{Thm3.1}(v),
using mainly (P7), was provided by Suyari \cite{Suyari2004} with a
correction by Ili\'{c} et al. \cite{IlicEtAl2013}. This generalization
considers the whole family of Tsallis entropies continuously depending on $q$
and includes the family of Havrda-Charv\'{a}t entropies \cite%
{HavrdaCharvat1967}. Other characterizations of the Tsallis entropy are due
to Abe \cite{Abe2004}, who uses a kind of conditional property, and to dos
Santos \cite{dosSantos1997}, who mainly refers to $q$-additivity.

Note that the statements listed in (i), (ii) and (iii) of Theorem \ref%
{Thm3.1} not covered by (iv) or (v) are well known and can easily be shown,
possibly with the exception of the statement in (ii) that the R\'{e}nyi
entropy $f_{q}$ does not satisfy (P7) nor (P8) for any $a>0$. To show this,
one can argue as follows. Given $q\neq 1$ and $a>0$, suppose that $f_{q}$
satisfies (P7) or (P8). Then%
\begin{equation*}
H\left( \frac{1}{4},\frac{1}{4},\frac{1}{4},\frac{1}{4}\right) =\left\{ 
\begin{array}{ll}
H\left( \frac{1}{2},\frac{1}{2}\right) +2\frac{1}{2^{a}}H\left( \frac{1}{2},%
\frac{1}{2}\right) & \text{by either (P7) and (P8),} \\ 
&  \\ 
2H\left( \frac{1}{2},\frac{1}{2}\right) & \text{by additivity,}%
\end{array}%
\right.
\end{equation*}%
which implies $a=1$; but $H=f_{q}$ does not satisfy (P7) nor (P8) for $a=1$.

The R\'{e}nyi entropy has been characterized in different way, however,
these characterizations are not as simple as those for the Tsallis entropy.
The first characterization, given by R\'{e}nyi \cite{Renyi1961}, uses mainly
additivity and a mean property for some special $q$-depending weighting, but
it is based on the extension of the system of distributions given by the set 
$\mathcal{P}$ to all incomplete distributions. Solving an open question
posed by R\'{e}nyi, Acz\'{e}l and Dar\'{o}czy \cite{AczelDaroczy1963} gave a
first characterization without considering incomplete distributions. Of
further characterizations, we only mention here one due to Jizba and
Arimitsu \cite{JizbaArimitsu2004}.

To conclude, let us underline the difficulty of finding the most general map 
$H:\mathcal{P}\rightarrow\lbrack0,\infty)$ satisfying a given set of
properties. To get a taste, the interested reader is referred to \cite%
{Hanel2011,Hanel2012}.

\subsection{Characterizations of the Kolmogorov-Sinai entropy}

In the same way as various Shannon entropy-like concepts have been
considered to quantify the diversity of probability distributions in
different situations, one can try to generalize the concept of
Kolmogorov-Sinai entropy by starting from entropies other than Shannon's.
The results of the last fifteen years show that this line of work does not
lead to anything new, at least for automorphisms.

In particular, attempts were made in \cite%
{HentschelProcaccia1983,GrassbergerProcaccia1984,EckmannRuelle1985} to
quantify the complexity of dynamical systems with the help of the R\'{e}nyi
entropy. Takens and Verbitskiy \cite{TakensVerbitskiy1998} have discussed
the R\'{e}nyi analogue of the Kolmogorov-Sinai entropy for variable $q>0$.
Given a measure-preserving dynamical system $(\Omega ,\mathcal{B},\mu ,T)$
and a finite partition $\alpha =\{A_{1},\ldots ,A_{n}\}$ of $\Omega $,
define the R\'{e}nyi entropy of $\alpha $ as%
\begin{equation*}
H_{q}(\alpha ):=\frac{1}{1-q}\ln \left( \sum_{i=1}^{n}\mu (A_{i})^{q}\right)
,
\end{equation*}%
the R\'{e}nyi entropy of $T$ with respect to $\alpha $ as%
\begin{equation*}
h_{\mu }^{q}(T,\alpha )=\liminf_{n\rightarrow \infty }\frac{1}{n}%
H_{q}(\alpha ^{(n)})
\end{equation*}%
(see Equation (\ref{alpha*(n)})), and, finally, the R\'{e}nyi entropy of $T$
as 
\begin{equation}
h_{\mu }^{q}(T)=\sup_{\alpha }h_{\mu }^{q}(T,\alpha ),  \label{RenyientT}
\end{equation}%
where the supremum is taken over all finite partitions $\alpha $ of $\Omega $%
. The surprising fact is that the new quantity (\ref{RenyientT}) is not
sensible in the case $q<1$ and coincides with the Kolmogorov-Sinai entropy
for $q\geq 1$.

\begin{theorem}
\cite{TakensVerbitskiy1998}\label{tvtheo} For each ergodic automorphism $T$
of a Lebesgue probability space $(\Omega ,\mathcal{B},\mu )$ with $h_{\mu
}(T)>0$, 
\begin{equation*}
h_{\mu }^{q}(T)=\left\{ 
\begin{array}{ll}
h_{\mu }(T) & \mbox{for }q\geq 1, \\ 
\infty & \mbox{else}.%
\end{array}%
\right.
\end{equation*}
\end{theorem}

Recall that a $\mu$-preserving map $T$ of a probability space $(\Omega ,%
\mathcal{B},\mu)$ is called automorphism if it is bijective and $T^{-1}$ is
also $\mu$-preserving.

Furthermore, Theorem \ref{tvtheo} remains true if $T$ is not an automorphism
or $h_{\mu}(T)=0$ \cite{TakensVerbitskiy1998}. On the other hand, as shown
in \cite{TakensVerbitskiy2002}, $h_{\mu}^{q}(T)$ can be strictly smaller
than $h_{\mu}(T)$ when ergodicity fails. Statements similar to Theorem \ref%
{tvtheo} for a Tsallis variant of the Kolmogorov-Sinai entropy were given by
Mes\'{o}n and Vericat \cite{MesonVericat2002}.

Another result in that direction found by Falniowski \cite{Falniowski2013}
sheds some light on the nature of the `ingredient' function $\eta
:[0,1]\rightarrow{\mathbb{R}}$ of the Kolmogorov-Sinai entropy given by $%
\eta(x)=-x\ln x$, $x>0$. Consider functions $g:[0,1]\rightarrow{\mathbb{R}}$
which are similar to $\eta$ in the sense that they are concave and satisfy 
\begin{equation}  \label{faln}
\lim_{x\rightarrow0^{+}}g(x)=g(0)=0\mbox{ for }x\in [0,1].
\end{equation}
For a partition $\alpha$ of $\Omega$, set 
\begin{equation*}
h_{\mu}(g,\alpha)=\sum_{A\in\alpha}g(\mu(A))
\end{equation*}
and 
\begin{equation*}
h_{\mu}(g,T,\alpha)=\limsup_{n\rightarrow\infty}\frac{1}{n}%
h_{\mu}(g,\alpha^{(n)}).
\end{equation*}
Then the \emph{$g$-entropy} of $T$, $h_{\mu}(g,T)$, is defined by 
\begin{equation*}
h_{\mu}(g,T)=\sup_{\alpha}h_{\mu}(g,T,\alpha).
\end{equation*}
It turns out that $g$-entropies with the same behavior near $0$ as $\eta$
provide the same information on the dynamical system $(\Omega,\mathcal{B}%
,\mu,T)$.

\begin{theorem}
\cite{Falniowski2013}\label{ftheo} Let $T$ be an ergodic automorphism of a
Lebesgue probability space $(\Omega,\mathcal{B},\mu)$ with $h_{\mu}(T)>0$,
and let $g:[0,1]\rightarrow{\mathbb{R}}$ be a concave function such that %
\eqref{faln} holds. If $c:=\limsup_{x\rightarrow0^{+}}\frac{g(x)}{\eta (x)}%
<\infty$, then 
\begin{equation*}
h_{\mu}(g,T)=c\,h_{\mu}(T),
\end{equation*}
where $h_{\mu}(T)$ is the Kolmogorov-Sinai entropy of $T$.
\end{theorem}

Additional information on the relation between the Kolmo\-go\-rov-Sinai
entropy and the $g$-entropies can be found in the paper of Falniowski. The
proof of Theorem \ref{ftheo} uses a very general characterization of the
Shannon entropy of finitely-valued stochastic processes given by Ornstein
and Weiss \cite{OrnsteinWeiss2012}.

A map $J$ from the class of finitely-valued ergodic stochastic processes $%
(X_{t})_{t\in {\mathbb{Z}}}$ into a metric space $(M,\rho )$ is said to be 
\emph{finitely observable} if there exist functions $S_{n}:{\mathbb{R}}%
^{n}\rightarrow M,\,n\in {\mathbb{N}}$ such that 
\begin{equation*}
\lim_{n\rightarrow \infty }S_{n}(X_{1}(\omega ),X_{2}(\omega ),\ldots
,X_{n}(\omega ))=J((X_{t})_{t\in {\mathbb{Z}}})\mbox{ for }\mu \mbox{-a.a. }%
\omega \in \Omega ,
\end{equation*}%
where $(X_{t})_{t\in {\mathbb{Z}}}$ is defined on a probability space $%
(\Omega ,\mathcal{B},\mu )$. $J$ is called an \emph{invariant} if it takes
the same values for finitely-valued stochastic processes with isomorphic
sequence space models (see Definition \ref{DefSupportX}).

\begin{theorem}
\label{Thm3.4}\cite{OrnsteinWeiss2012} Every finitely observable invariant
for the class of all finitely-valued ergodic processes is a continuous
function of the Kolmogorov-Sinai entropy.
\end{theorem}

A consequence of this statement is, roughly speaking, that the only useful
complexity measures for dynamical systems are continuous functions of the
Kolmogorov-Sinai entropy.

\section{Other entropies}

Of the various entropies and entropy-like quantities proposed in the
literature, the selection below responds to both their applications and
mathematical insights.

\subsection{Approximate and sample entropy}

Approximate entropy and sample entropy are widely-used quantities for
measuring complexity of (finite) time series $(x_{t})_{t=1}^{N}$ and the
underlying dynamical systems. They quantify the change in the relative
frequencies of length $k$ time-delay vectors 
\begin{equation*}
(x_{t})_{t=j}^{j+(k-1)\tau},\;\; j=1,2,\ldots,N-(k-1)\tau,
\end{equation*}
with increasing $k$. Here, for simplicity, we assume delay time $\tau=1$,
but searching for delay times with an optimal perception of the dynamical
structure, is an important task in nonlinear data analysis; see \cite[Sect.
9.2]{Kantz2005} for several methods.

We first recall the concepts of approximate entropy and sample entropy from 
\cite{Pincus1991,RichmanMoorman2000} and then discuss where they come from
(see \cite{UnakafovaThesis} for a more extended discussion).

\begin{definition}
\label{ApEnDef} Given a time series $(x_{t})_{t=1}^{N}$, and a tolerance $%
\epsilon >0$ for accepting time-delay vectors of length $k\in \mathbb{N}$ as
similar, the \emph{approximate entropy} is defined as 
\begin{equation}
\mathrm{ApEn}(k,\epsilon ,(x_{t})_{t=1}^{N})=\Phi (k,\epsilon
,(x_{t})_{t=1}^{N})-\Phi (k+1,\epsilon ,(x_{t})_{t=1}^{N}),  \label{ApEn}
\end{equation}%
where 
\begin{equation}
\Phi (k,\epsilon ,(x_{t})_{t=1}^{N})=\frac{1}{N-k+1}\sum_{i=1}^{N-k+1}\ln (%
\widehat{C}(i,k,\epsilon ,(x_{t})_{t=1}^{N}))  \label{Phi}
\end{equation}%
and%
\begin{equation}
\widehat{C}(i,k,\epsilon ,(x_{t})_{t=1}^{N})=\frac{\left\vert \{1\leq j\leq
N-k+1:\max\limits_{l=0,1,\ldots ,k-1}{\lvert x_{i+l}-x_{j+l}\rvert }\leq
\epsilon \}\right\vert }{N-k+1},  \label{cform}
\end{equation}%
and the \emph{sample entropy} is defined as 
\begin{equation}
\mathrm{SampEn}(k,\epsilon ,(x_{t})_{t=1}^{N})=\ln \widehat{C}(k,\epsilon
,(x_{t})_{t=1}^{N})-\ln \widehat{C}(k+1,\epsilon ,(x_{t})_{t=1}^{N}),
\label{SampEn}
\end{equation}%
where 
\begin{equation}
\widehat{C}(k,\epsilon ,(x_{t})_{t=1}^{N})=\frac{2\left\vert \{(i,j):1\leq
i<j\leq N-k+1,\max\limits_{l=0,1,\ldots ,k-1}{\lvert x_{i+l}-x_{j+l}\rvert }%
\leq \epsilon \}\right\vert }{(N-k-1)(N-k)}.  \label{c2form}
\end{equation}
\end{definition}

Thus, with respect to the maximum norm, $\widehat{C}(i,k,\epsilon
,(x_{t})_{t=1}^{N})$ is the relative frequency of the vectors $%
(x_{t})_{t=j}^{j+m-1}$, $1\leq j\leq N-k+1$, which are within a distance $%
\epsilon $ from the vector $(x_{t})_{t=i}^{i+k-1}$, and $\widehat{C}%
(k,\epsilon ,(x_{t})_{t=1}^{N})$ is the relative frequency of pairs of
vectors $(x_{t})_{t=i}^{i+m-1},(x_{t})_{t=j}^{j+m-1}$, $1\leq i,j\leq N-k+1$%
, within a distance $\epsilon $. At a first glance, the quantities
considered do not look much like entropies. In the following, we justify (i)
the approximate entropy as an estimate of the Kolmogorov-Sinai entropy,
referring to a discussion in \cite{EckmannRuelle1985} based on ideas from 
\cite{GrassbergerProcaccia1983a,GrassbergerProcaccia1983}, and (ii) the
sample entropy as an estimate of the $H_{2}(T)$-entropy, introduced,
according to \cite{BroerTakens2011}, in \cite{Takens1981}.

The correlation integral, which is necessary to understand what is really
behind approximate and sample entropies, was originally defined in \cite%
{GrassbergerProcaccia1983a}. Given a probability space, one can interpret
the correlation integral as the probability of two points being closer than
a distance $\epsilon$.

\begin{definition}
Let $(\Omega ,\rho )$ be a metric space with a Borel probability measure $%
\mu $, and let $\epsilon >0$. Then the \textit{correlation integral} is
given by 
\begin{equation}
C(\epsilon )=(\mu \times \mu )\{(\omega ,\omega ^{\prime })\in \Omega \times
\Omega :\rho (\omega ,\omega ^{\prime })<\epsilon \}.
\end{equation}
\end{definition}

This definition is equivalent to 
\begin{equation*}
C(\epsilon )=\int_{\Omega }\mu (B(\omega ,\epsilon ))\mathrm{d}\mu (\omega )%
\text{,}
\end{equation*}%
where $B(\omega ,\epsilon )$ is the ball of center $\omega $ and radius $%
\epsilon $ in $(\Omega ,\rho )$, i.e.,%
\begin{equation*}
B(\omega ,\epsilon )=\{\omega ^{\prime }\in \Omega :\rho (\omega ,\omega
^{\prime })<\epsilon \}.
\end{equation*}%
Indeed, if 
\begin{equation*}
1_{\delta _{\epsilon }}(\omega ,\omega ^{\prime })=%
\begin{cases}
1 & \text{if }\rho (\omega ,\omega ^{\prime })<\epsilon , \\ 
0 & \text{otherwise,}%
\end{cases}%
\end{equation*}%
then, by Fubini's theorem \cite{Choe2005}, 
\begin{equation*}
\int_{\Omega \times \Omega }1_{\delta _{\epsilon }}(\omega ,\omega ^{\prime
})\mathrm{d}(\mu \times \mu )=\int_{\Omega }\left( \int_{\Omega }1_{\delta
_{\epsilon }}(\omega ,\omega ^{\prime })\mathrm{d}\mu (\omega ^{\prime
})\right) \mathrm{d}\mu (\omega )=\int_{\Omega }\mu (B(\omega ,\epsilon ))%
\mathrm{d}\mu (\omega ).
\end{equation*}

Let us consider a Borel-measurable map $T:\Omega \rightarrow \Omega $.
Furthermore, let $B_{k}(\omega ,\epsilon )$ be the ball of center $\omega $
and radius $\epsilon $ with respect to the metric $\rho _{k}$ defined in
Equation (\ref{rhok}). The \emph{correlation integral} $C(k,\epsilon )$ over
the balls $B_{k}(\omega ,\epsilon )$ in $\Omega $ is then given by 
\begin{equation}
C(k,\epsilon )=\int_{\Omega }\mu (B_{k}(\omega ,\epsilon ))\mathrm{d}\mu
(\omega ).  \label{MetricDefinition2}
\end{equation}

The following deep result of Brin and Katok \cite{BrinKatok1983} is the base
for relating the approximate and sample entropies to the Komogorov-Sinai
entropy.

\begin{theorem}
\label{BrinKatok} Let $(\Omega ,B,\mu ,T)$ be a measure-preserving dynamical
system, where $\Omega $ is a compact metric space, $\mu $ is a non-atomic
Borel probability measure, and $T$ is a continuous map. Then 
\begin{equation}
\int_{\Omega }\lim_{\epsilon \rightarrow 0}\underset{k\rightarrow \infty }{%
\lim \sup }-\frac{1}{k}\ln \mu (B_{k}(\omega ,\epsilon ))\mathrm{d}\mu
(\omega )=\int_{\Omega }\lim_{\epsilon \rightarrow 0}\underset{k\rightarrow
\infty }{\lim \inf }-\frac{1}{k}\ln \mu (B_{k}(\omega ,\epsilon ))\mathrm{d}%
\mu (\omega )=h_{\mu }(T).  \label{Eq_Brin_Katok}
\end{equation}
\end{theorem}

Recall that a measure $\mu$ is called \textit{non-atomic} if for any
measurable set $A$ with $\mu(A)>0$ there exists a measurable set $B\subset A$
such that $\mu(A)>\mu(B)>0$.

Further, the following statement can be shown (see \cite{UnakafovaThesis}).

\begin{theorem}
\label{MyTheorem} Under the assumptions of Theorem~\ref{BrinKatok} for an
ergodic $T$, if for all $\epsilon >0$ the limit 
\begin{equation}
\lim_{k\rightarrow \infty }\lim_{n\rightarrow \infty }\frac{1}{n}\left(
\sum_{j=0}^{n-1}\ln \mu (B_{k}(T^{j}(\omega ),\epsilon
))-\sum_{j=0}^{n-1}\ln \mu (B_{k+1}(T^{j}(\omega ),\epsilon ))\right)
\label{LimitExistence}
\end{equation}%
exists, then it holds 
\begin{equation}
h_{\mu }(T)= \lim_{\epsilon \rightarrow 0}\lim_{k\rightarrow \infty
}\lim_{n\rightarrow \infty }\frac{1}{n}\left( \sum_{j=0}^{n-1}\ln \mu
(B_{k}(T^{j}(\omega ),\epsilon ))-\sum_{j=0}^{n-1}\ln \mu
(B_{k+1}(T^{j}(\omega ),\epsilon ))\right) .  \label{ApEn2}
\end{equation}
\end{theorem}

\noindent In general it is not known whether the limit \eqref{LimitExistence}
exists.

In this setting, assume that $x_{t}=X(T^{t}(\omega ))$, where $X:\Omega
\rightarrow \mathbb{R}$ is a Borel-measurable function, i.e., the time
series $(x_{t})_{t=1}^{N}$ is obtained from the dynamical system $(\Omega ,%
\mathcal{B},\mu ,T)$ via an observable $X$. Here $\mathcal{B}$ is the Borel $%
\sigma $-algebra of $\Omega $, and $\mu $ is a Borel measure. Then, for an
ergodic $T$ the term \eqref{cform} is a natural estimate of $\mu
(B_{k}(T^{i}(\omega ),\epsilon ))$. This being the case, comparison of the
definition (\ref{ApEn}) of approximate entropy with the right hand side of (%
\ref{ApEn2}) for finite $k$, $n$ and $\epsilon $, leads to the conclusion
that $\mathrm{ApEn}(k,\epsilon ,(x_{t})_{t=1}^{N})$ approximates the
Kolmogorov-Sinai entropy under the conditions of Theorem \ref{MyTheorem}.

The switch from the metric $\rho $ of the state space $\Omega $ to the
maximum norm metric $\rho _{k}$ of the time delay vectors $%
(x_{i},x_{i+1},...,x_{i+k-1})$, Equation (\ref{rhok}), is justified by
Takens' theorem in the case of dynamical systems with certain natural
differentiability assumptions. This theorem guarantees that generically the
space of delay vectors `reconstructs' the dynamical system with equivalence
of the original metric and the maximum norm metric (see \cite%
{BroerTakens2011}, in particular Chapter 6, for more details).

Next we consider the relationship between the $H_{2}(T)$-entropy (see below)
and its estimation by means of the correlation integral.

\begin{definition}
\label{CorrelationEntropy} Let $(\Omega ,\mathcal{B},\mu ,T)$ be a
measure-preserving dynamical system, where $\Omega $ is a compact metric
space, $\mu $ is a Borel probability measure, and $T$ is a continuous map.
The \textit{$H_{2}(T)$-entropy} and the \textit{correlation entropy $\mathrm{%
CE}(T,1)$} are then defined by 
\begin{align}
H_{2}(T)& =\underset{\epsilon \rightarrow 0}{\lim \sup \,\,}\underset{%
k\rightarrow \infty }{\lim \sup \,\,}-\frac{1}{k}\ln \int_{\Omega }\mu
(B_{k}(\omega ,\epsilon ))\mathrm{d}\mu (\omega ), \\
\mathrm{CE}(T,1)& =\lim_{\epsilon \rightarrow 0}\lim_{k\rightarrow \infty }-%
\frac{1}{k}\int_{\Omega }\ln \mu (B_{k}(\omega ,\epsilon ))\mathrm{d}\mu
(\omega ).  \label{CEdef}
\end{align}
\end{definition}

Note that the limit for $\epsilon \rightarrow 0$ in \eqref{CEdef} exists due
to the monotonicity properties of $-\frac{1}{k}\int_{\Omega }\ln \mu
(B_{k}(\omega ,\epsilon ))\mathrm{d}\mu (\omega )$; see \cite[Lemma~2.1]%
{Verbitskiy2000} for details. The limit for $k\rightarrow \infty $ in %
\eqref{CEdef} exists by Lemma 2.14 from \cite{Verbitskiy2000}.

The $H_{2}(T)$-entropy was introduced in \cite{Takens1981} (see also \cite%
{BroerTakens2011}). Inspired by \cite%
{EckmannRuelle1985,GrassbergerProcaccia1984,Takens1981}, where the R\'{e}nyi
approach was used, Takens and Verbitskiy \cite{TakensVerbitskiy1998}
introduced a correlation entropy $\mathrm{CE}(T,q)$ of order $q$. Here we
will only consider the case $q=1$. For an interesting study of properties of
the correlation entropies $\mathrm{CE}(T,q)$, we refer also to Verbitskiy's
dissertation \cite{Verbitskiy2000}.

The following theorem relates $H_{2}(T)$-entropy, correlation entropy and
Kolmogo\-rov-Sinai entropy (for a proof based on ideas and results from \cite%
{TakensVerbitskiy1998,Verbitskiy2000,BroerTakens2011}, see \cite%
{UnakafovaThesis}).

\begin{theorem}
\label{CombinedCorrelation} Under the conditions of Theorem~\ref{BrinKatok}
it holds 
\begin{equation}  \label{ComCorr}
H_{2}(T)\leq\mathrm{CE}(T,1)=h_{\mu}(T).
\end{equation}
\end{theorem}

\noindent Note that the inequality in \eqref{ComCorr} follows directly from
Jensen's inequality.

Applying Birkhoff's ergodic theorem twice, one obtains 
\begin{eqnarray*}
C(k,r)\hspace*{-8mm} &=&\hspace*{-8mm}\int_{\Omega }\mu (B_{k}(\omega ,r))\,%
\mathrm{d}\mu (\omega ) \\
&\underset{\text{ for a.a. }\omega _{1}\in \Omega }{=}&\lim_{n\rightarrow
\infty }\frac{1}{n}\sum_{j=0}^{n-1}\mu (B_{k}(T^{j}(\omega _{1}),r)) \\
&\underset{\text{ for a.a. }\omega _{2}\in \Omega }{=}&\lim_{n\rightarrow
\infty }\frac{1}{n^{2}}\sum_{i=0}^{n-1}\sum_{j=0}^{n-1}1_{\delta
_{r,k}}(T^{i}(\omega _{1}),T^{j}(\omega _{2})),
\end{eqnarray*}%
where 
\begin{equation*}
1_{\delta _{\epsilon ,k}}(T^{i}(\omega _{1}),T^{j}(\omega _{2}))=%
\begin{cases}
1, & \text{if }\rho _{k}(T^{i}(\omega _{1}),T^{j}(\omega _{2}))<\epsilon ,
\\ 
0, & \text{otherwise}.%
\end{cases}%
\end{equation*}

Equation \eqref{c2form} turns out to be an appropriate estimate of the
correlation integral from a time series $(x_{i})_{i\in \mathbb{N}%
}=(X(T^{i}(\omega )))_{i\in \mathbb{N}}$, again under consideration of
Takens' reconstruction theorem (see \cite{Takens1981,
GrassbergerProcaccia1983,BroerTakens2011}). In \cite{DenkerKeller1986} it
was shown that for almost all $\omega \in \Omega $, 
\begin{equation*}
\lim_{N\rightarrow \infty }\widehat{C}(k,\epsilon ,(X(T^{t}(\omega
)))_{t=1}^{N})=C(k,\epsilon ).
\end{equation*}%
We refer to \cite{BroerTakens2011,Borovkova1998} for a further discussion on
the estimation of the correlation integral.

Finally, the following estimate of $H_{2}(T)$-entropy was proposed in \cite%
{Takens1981}, and then applied in \cite%
{GrassbergerProcaccia1983,GrassbergerProcaccia1983a} (see also \cite%
{BroerTakens2011}): 
\begin{equation}
\widehat{H}_{2}(k,\epsilon ,(x_{t})_{t=1}^{N})=\ln \frac{\widehat{C}%
(k,\epsilon ,(x_{t})_{t=1}^{N})}{\widehat{C}(k+1,\epsilon ,(x_{t})_{t=1}^{N})%
}.  \label{H2estimate}
\end{equation}%
Equation \eqref{H2estimate} is nothing else but the definition \eqref{SampEn}
of sample entropy.

\subsection{Directional entropy}

Directional entropy was introduced by Milnor in \cite{Milnor1988} for
cellular automaton maps. In its simplest formulation, a one-dimensional
cellular automaton map is a transformation $F:S^{\mathbb{Z}}\rightarrow S^{%
\mathbb{Z}}$, where $S$ is a finite set of symbols, such that $(F(%
\boldsymbol{s}))_{i}$, $i\in\mathbb{Z}$, depends only locally on the
components of $\mathbf{s}=(s_{i})_{i\in\mathbb{Z}}$. Specifically, there is
an $r\geq1$, called the radius of $F$, and a `local rule' $%
f:S^{2r+1}\rightarrow S$ such that 
\begin{equation}
(F(\boldsymbol{s}))_{i}=f(s_{i-r},...,s_{i+r}).  \label{directional1}
\end{equation}
The transformation $F$ implements the state update of the cellular automaton
(CA) $(S^{\mathbb{Z}}, F)$, so one writes $\mathbf{s}_{t}=F^{t}(\mathbf{s})$%
. Let us mention in passing that a map between two shift spaces of the form (%
\ref{directional1}) is called a block map. They are characterized by being
continuous and shift commuting \cite{Hedlund1969}.

This situation generalizes readily to $d$-dimensional CA; here $F:S^{\mathbb{%
Z}^{d}}\rightarrow S^{\mathbb{Z}^{d}}$. Cellular automata can be considered
as toy models of spatially extended systems ($d$-dimensional space being
modelled by the lattice $\mathbb{Z}^{d}$) because the time evolution of each
component $s_{i\in\mathbb{Z}^{d}}$ (or `cell') depends on the states of the
neighboring components. Roughly speaking, Milnor's topological and
measure-theoretic directional entropy measure the complexity of such a
spatiotemporal dynamics in the direction of a vector $\mathbf{v}\in\mathbb{R}%
^{d}$.

Following \cite{Robinson2011}, let $V$ be an $k$-dimensional vector subspace
of $\mathbb{R}^{d}$, $1\leq k<d$, and let $Q\subset V$ be a unit cube in $V$
built on an orthonormal basis of $V$. Likewise, let $Q^{\prime }\subset
V^{\bot }$ be such a unit cube in the orthogonal subspace of $V$, centered
at the origin. The set $W_{a,b}:=(aQ+bQ^{\prime })\cap \mathbb{Z}^{d}$, $%
a,b\in \mathbb{R}$, is called a \textit{window}. Furthermore, let $T$ be a
measure-preserving $\mathbb{Z}^{d}$ action on a Lebesgue probability space $%
(\Omega ,\mathcal{B},\mu )$, and $\alpha $ a finite partition of $\Omega $.
Define,%
\begin{equation*}
\alpha _{V,a,b}=\bigvee\limits_{\mathbf{n}\in W_{a,b}}T^{-\mathbf{n}}\alpha ,
\end{equation*}

\begin{equation}
h_{\mu,k}(T,V,\alpha,b)=\lim_{a\rightarrow\infty}\frac{1}{a^{k}}%
H_{\mu}(\alpha_{V,a,b}),  \label{directional2}
\end{equation}
and%
\begin{equation}
h_{\mu,k}(T,V,\alpha)=\lim_{b\rightarrow\infty}h_{\mu,k}(T,V,\alpha,b).
\label{directional3}
\end{equation}
The convergence of the limits in (\ref{directional2}) and (\ref{directional3}%
) was proved by Milnor in \cite{Milnor1988}.

\begin{definition}
If $1\leq k<d$, the $k$-\emph{dimensional directional entropy in the
direction} $V$ is 
\begin{equation}
h_{\mu,k}(T,V)=\sup_{\alpha}h_{\mu,k}(T,V,\alpha).  \label{directent}
\end{equation}
\end{definition}

The directional entropy (\ref{directent}) does not depend on the ortonormal
bases used to define the windows $W_{a,b}$. If $k=d$ (i.e., $V=\mathbb{R}%
^{d} $), then the windows have the form $W_{a}:=aQ$, and $h_{\mu ,d}(T,%
\mathbb{R}^{d})=h_{\mu ,d}(T)$, Equation (\ref{action3}). As for the
relation between different directional entropies of the same action, the
following is an interesting result. If $V_{1}\subset V_{2}$ are subspaces of 
$\mathbb{R}^{d}$ with $k_{1}=\dim V_{1}<k_{2}=\dim V_{2}\leq d$, and $h_{\mu
,k_{2}}(T,V_{2})>0$, then $h_{\mu ,k_{1}}(T,V_{1})=\infty $ \cite%
{Robinson2011}.

For the directional entropy, including $h_{\mu,d}(T)$, a similar result to
Theorem \ref{Thm2.8} holds.

\begin{theorem}
Let $1\leq k\leq d$. If $(\alpha _{n})_{n\in \mathbb{N}}\nearrow \varepsilon 
$ (the partition of $\Omega $ into separate points), then%
\begin{equation}
h_{\mu ,k}(T,V)=\lim_{n\rightarrow \infty }h_{\mu ,k}(T,V,\alpha _{n}).
\label{directional5}
\end{equation}
\end{theorem}

The interested reader is referred to \cite%
{Boyle1997,Courbage2006,Afraimovich2015,Miles2015} for further insights into
directional entropy.

\subsection{Permutation entropy}

Let $I$ be a closed interval of $\mathbb{R}$ endowed with the Borel $\sigma $%
-algebra $\mathcal{B}$. A map $T:I\rightarrow I$ is said to be \textit{%
piecewise}, \textit{strictly monotone} if there is a finite partition of $I$
into intervals such that, on each of those intervals, $T$ is continuous and
strictly monotone. Furthermore, let $r:\{0,1,...,L-1\}\rightarrow
\{0,1,...,L-1\}$ be the permutation sending $i$ to $r_{i}$, $0\leq i\leq L-1$%
. For notational convenience, we will write $r=(r_{0},r_{1},...,r_{L-1})$.

\begin{definition}
\label{DefA}We say that $x\in I$ has \emph{ordinal }$L$\emph{-pattern} $r$,
or that $x$ is of \emph{type} $r$, if%
\begin{equation}
T^{r_{0}}(x)<T^{r_{1}}(x)<...<T^{r_{L-1}}(x).  \label{Pr}
\end{equation}
\end{definition}

We refer to $L$ as the length of the ordinal pattern $r$. Other conventions,
such as reversing the order relations in (\ref{Pr}), can be also found in
the literature.

Let us highlight at this point that the permutations of $L$ elements (or
ordinal $L$-patterns for that matter) build the \textit{symmetric group of
degree} $L$, denoted by $\mathcal{S}_{L}$, when endowed with the `product'
(actually, composition of permutations)%
\begin{equation}
r\circ s =(r_{0},r_{1},...,r_{L-1})\circ(s_{0},s_{1},...,s_{L-1})
=(s_{r_{0}},s_{r_{1}},...,s_{r_{L-1}}),  \label{ros}
\end{equation}
for all $r,s\in\mathcal{S}_{L}$, the unity being the identity permutation $%
(0,1,...,L-1)$. $\mathcal{S}_{L}$ is non-commutative for $L\geq3$, and its
cardinality is $L!$. The algebraic structure of ordinal patterns will be
exploited in the next section.

Let $P_{r}$ be the set of $x\in I$ of type $r\in \mathcal{S}_{L}$, i.e.,
points $x$ satisfying (\ref{Pr}), and%
\begin{equation}
\pi _{L}=\{P_{r}\neq \emptyset :r\in \mathcal{S}_{L}\}.  \label{pi_L}
\end{equation}%
That is, $\pi _{L}$ is a partition of $I$ into subsets whose points have the
same ordinal $L$-pattern, possibly except for a set of Lebesgue measure 0
(the set of periodic points with period less than $L$). Alternatively, one
can agree some convention to rank $T^{r_{i}}(x)$ and $T^{r_{j}}(x)$, $0\leq
r_{i},r_{j}\leq L-1$, whenever $T^{r_{i}}(x)=T^{r_{j}}(x)$. $\pi _{L}$ is
called an \textit{ordinal partition}. Note that $\pi _{L+1}$ is a refinement
of $\pi _{L}$.

\begin{definition}
\label{DefB} Let $T:I\rightarrow I$ be a piecewise, strictly monotone map,
and $\mu$ be a $T$-invariant measure on $(I,\mathcal{B})$. Then, the \emph{%
metric permutation entropy} of $T$ is defined as%
\begin{equation}
h_{\mu}^{\ast}(T)=-\lim_{L\rightarrow\infty}\frac{1}{L}H_{\mu}(\pi_{L})=-%
\lim_{L\rightarrow\infty}\frac{1}{L}\sum_{P_{r}\in\pi_{L}}\mu(P_{r})\ln
\mu(P_{r}).  \label{h*(T)}
\end{equation}
Furthermore, the \emph{topological permutation entropy} of $T$ is defined as%
\begin{equation}
h_{0}^{\ast}(T)=\lim_{L\rightarrow\infty}\frac{1}{L}\ln\left\vert \pi
_{L}\right\vert .  \label{h*0(T)}
\end{equation}
\end{definition}

Note that the definitions (\ref{h*(T)}) and (\ref{h*0(T)}) are simpler than
their standard counterparts (\ref{h(T)}) and (\ref{h_top(T)}), respectively,
since no supremum must be taken. For this reason the following result (also
showing the convergence of (\ref{h*(T)}) and (\ref{h*0(T)})) is remarkable.

\begin{theorem}
\cite{Bandt2002A} \label{ThmC} It holds: (i) $h_{\mu}^{\ast}(T)=h_{\mu}(T)$,
and (ii) $h_{0}^{\ast}(T)=h_{top}(T).$
\end{theorem}

Metric permutation entropy and, especially, the \textit{metric permutation
entropy of order $L$},%
\begin{equation*}
h_{\mu}^{\ast}(T,L)=H_{\mu}(\pi_{L})=-\sum_{P_{r}\in\pi_{L}}\mu(P_{r})\ln%
\mu(P_{r}),
\end{equation*}
have a number of interesting applications, in particular to time series
analysis (complexity, change of dynamical regime, segmentation,
discrimination,...). Some of them will be reviewed in Sect. 5.1. In this
regard, the \textit{conditional entropy of ordinal patterns},%
\begin{equation*}
h_{\mu,cond}^{\ast}(T,L)=H_{\mu}(T^{-1}\pi_{L}\left\vert \pi_{L}\right)
:=H_{\mu}(\pi_{L}\vee T^{-1}\pi_{L})-H_{\mu}(\pi_{L}),
\end{equation*}
has proved to be useful, also in the estimation of the metric entropy \cite%
{Unakafov2014}.

As for the \textit{topological permutation entropy of order $L$},%
\begin{equation*}
h_{0}^{\ast }(T,L)=\frac{1}{L}\ln \left\vert \pi _{L}\right\vert ,
\end{equation*}%
it amounts to counting ordinal $L$-patterns. Application examples to time
series analysis can be found in \cite{Amigo2015}. Moreover, Theorem \ref%
{ThmC}\textit{(ii)}, i.e.,%
\begin{equation*}
\left\vert \pi _{L}\right\vert \propto e^{h_{top}(T)L},
\end{equation*}%
along with $\left\vert S_{L}\right\vert =L!\propto e^{L\ln L}$, has an
interesting offshoot. Namely, the existence of \textit{forbidden ordinal
patterns}, that is, ordinal patterns that cannot be realized by any orbit of
the system. For example, the ordinal pattern $(2,1,0)$ cannot be realized by
the logistic map $f(x)=4x(1-x)$, $0\leq x\leq 1$, i.e., there is no initial
condition $x\in \lbrack 0,1]$ such that $f^{2}(x)<f(x)<x$ \cite%
{Amigo2006,Amigo2008A}. It is trivial that if $L_{\min }$ is the minimum
length of the forbidden patterns, then there are forbidden $L$-patterns for
all $L\geq L_{\min }$ (in fact, a superexponentially growing number of
them). This result was used in \cite{Amigo2008B,Amigo2010A,Amigo2010B} to
discriminate deterministic noisy signals from white noise. We come back to
this application in Sect. 5.4.

The computation of permutation entropies of finite orders benefits from the
practical advantages of using ordinal patterns. Thus, $h_{\mu}^{\ast}(T,L)$
and $h_{0}^{\ast}(T,L)$ are computationally fast and relatively robust
against observational noise. Actually, ordinal patterns of deterministic
signals are more robust than those of random signals due to a mechanism
called dynamic robustness \cite{Amigo2010B}. Furthermore, the calculation of
ordinal patterns does not require a prior knowledge of the data range, what
allows real-time data processing and analysis.

Metric permutation entropy was extended to general dynamical systems along
two different lines \cite{AmigoKeller2013}. One of them follows the
standard, partition-based approach to metric entropy and turns out to be
equivalent to it \cite{Amigo2012C}. The other generalization considers
real-valued maps, `observables' on dynamical sets, and ordinal patterns
defined on their domain. It can be shown that under some separation
conditions of the maps, this approach leads also to the metric entropy \cite%
{Keller2012, AntonioukEtAl2014, KellerEtAl2015}.

An open problem is the extension of topological permutation entropy to more
general settings, as successfully done with the metric version. Misiurewicz
proved in \cite{Misiurewicz2003} that Theorem \ref{ThmC}\textit{(ii)} does
not hold if $T$ is not piecewise, strictly monotone, so no straightforward
generalization can be expected. Incidentally, Theorem \ref{ThmC}\textit{(ii)}
shows that, under some restrictions, topological entropy can be computed via
partitions instead of open covers.

\bigskip

\subsection{Transfer entropy}

A relevant question in time series analysis of coupled random or
deterministic processes is the causality relation, i.e., which process is
driving and which is responding. A possible method to discriminate cause
from effect, proposed originally by Wiener \cite{Wiener1956}, is as follows:
\textquotedblleft For two simultaneously measured signals, if we can predict
the first signal better by using the past information from the second one
than by using the information without it, then we call the second signal
causal to the first one\textquotedblright. The implementation of this
principle in linear time series analysis via autoregressive processes goes
by the name of Granger causality \cite{Granger1969}. Transfer entropy,
introduced by Schreiber in \cite{Schreiber2000} for measuring the
information exchanged between two processes in both directions separately,
can be considered an information-theoretic implementation of Wiener
causality.

Let $\mathbf{X}=(X_{t})_{t\in \mathbb{T}}$ and $\mathbf{Y}=(Y_{t})_{t\in 
\mathbb{T}}$ be two stationary random processes. Then, the transfer entropy
from $\mathbf{Y}$ to $\mathbf{X}$, $T_{\mathbf{Y}\rightarrow \mathbf{X}}$,
is the reduction of uncertainty in future values of $\mathbf{X}$, given past
values of $\mathbf{X}$, due to the additional knowledge of past values of $%
\mathbf{Y}$. Set $X_{t}^{(n)}:=X_{t},...,X_{t-n+1}$, and $%
Y_{t}^{(k)}:=Y_{t},...,Y_{t-k+1}$. To be more specific:

\begin{definition}
\label{Definition4.12}\cite{Schreiber2000} The \emph{transfer entropy from
the process }$\mathbf{Y}$\emph{\ to the process }$\mathbf{X}$ with coupling
delay $\Lambda \geq 1$ is defined as 
\begin{align}
T_{\mathbf{Y}\rightarrow \mathbf{X}}(\Lambda )& =H\left( X_{t+\Lambda
}\right\vert X_{t}^{(n)})-H\left( X_{t+\Lambda }\right\vert
X_{t}^{(n)},Y_{t}^{(k)})  \label{Transfer1} \\
& =\sum_{x_{t+\Lambda },x_{t}^{(n)},y_{t}^{(k)}}p(x_{t+\Lambda
},x_{t}^{(n)},y_{t}^{(k)})\ln \frac{p\left( x_{t+\Lambda }\right\vert
x_{t}^{(n)},y_{t}^{(k)})}{p\left( x_{t+\Lambda }\right\vert x_{t}^{(n)})}. 
\notag
\end{align}
\end{definition}

Therefore, $T_{\mathbf{Y}\rightarrow\mathbf{X}}(\Lambda)=0$ if and only if $%
p\left( x_{t+\Lambda}\right\vert x_{t}^{(n)},y_{t}^{(k)})=p\left(
x_{t+\Lambda}\right\vert x_{t}^{(n)})$ for all $x_{t+%
\Lambda},x_{t}^{(n)},y_{t}^{(k)}$. In physical terms, the transfer entropy
from $\mathbf{Y}$ to $\mathbf{X}$ vanishes if and only if $\mathbf{Y}$ has
no influence on $\mathbf{X}$ or both processes are fully synchronized.
Otherwise, $T_{\mathbf{Y}\rightarrow\mathbf{X}}(\Lambda)>0$ and we say that
there is an information transfer from the process $\mathbf{Y}$ to the
process $\mathbf{X}$. In \cite{Schreiber2000} and often in applications, $%
\Lambda=1$.

For notational simplicity, a possible dependence of $T_{\mathbf{Y}%
\rightarrow \mathbf{X}}$ on the history lengths $n$ and $k$ is not
explicitly shown. Note from the definition of conditional entropy, Equation (%
\ref{Hcond}), and the definition of mutual information, Equation (\ref%
{I(X;Y)}), that $T_{\mathbf{Y}\rightarrow \mathbf{X}}(\Lambda )$ is actually
a conditional mutual information, to wit:%
\begin{align}
T_{\mathbf{Y}\rightarrow \mathbf{X}}(\Lambda )& =I(X_{t+\Lambda
};Y_{t}^{(k)}\mid X_{t}^{(n)})  \label{Transfer2} \\
& =\sum_{x_{t+\Lambda },x_{t}^{(n)},y_{t}^{(k)}}p(x_{t+\Lambda
},x_{t}^{(n)},y_{t}^{(k)})\ln \frac{p(x_{t+\Lambda },y_{t}^{(k)}\mid
x_{t}^{(n)})}{p(x_{t+\Lambda }\mid x_{t}^{(n)})\,p(y_{t}^{(k)}\mid
x_{t}^{(n)})}.  \notag
\end{align}

Unlike the (unconditioned) mutual information (Equation (\ref{I(X;Y)B})), $%
T_{\mathbf{Y}\rightarrow \mathbf{X}}(\Lambda )$ is asymmetric under the
exchange of the processes $\mathbf{X}$ and $\mathbf{Y}$. Therefore, the 
\textit{directionality indicator} 
\begin{equation}
\Delta T_{\mathbf{Y}\rightarrow \mathbf{X}}(\Lambda )=T_{\mathbf{Y}%
\rightarrow \mathbf{X}}(\Lambda )-T_{\mathbf{X}\rightarrow \mathbf{Y}%
}(\Lambda )=-\Delta T_{\mathbf{X}\rightarrow \mathbf{Y}}(\Lambda ).
\label{DeltaTE}
\end{equation}%
measures the \textit{net} transfer of information between the processes $%
\mathbf{X}$ and $\mathbf{Y}$. For example, if $\Delta T_{\mathbf{Y}%
\rightarrow \mathbf{X}}(\Lambda )>0$ then $\mathbf{Y}$ is the driving
process with a coupling delay $\Lambda $. This so-called coupling direction
is one of the main objectives in the study of interacting systems.

Other definitions of transfer entropy (e.g., with $t-1$ instead of $t$, or
with $k=n$) can be also found in the literature. For a more general notion,
called \textit{momentary information transfer}, see \cite{Pompe2011}.

In practice, data are finite if not sparse. In the latter case, one often
uses transfer entropy with the least dimension,%
\begin{equation}
T_{\mathbf{Y}\rightarrow \mathbf{X}}(\Lambda )=I(X_{t+\Lambda };Y_{t}\mid
X_{t})  \label{Transfer3}
\end{equation}%
(i.e., $k=n=1$), as we assume hereafter for simplicity. Furthermore, to
estimate transfer entropy (as well as other statistics and observables) one
resorts to symbolic representations. Indeed, this is a common technique in
time series analysis which consists in trading-off realizations of a random
variable for symbols belonging to some convenient alphabet. Instances of
this procedure are binning in traditional statistical data analysis,
instantaneous phases via the Hilbert transform, and symbolic dynamics in
nonlinear time series analysis (Sect. 2.3). The rationale of the latter
instance is that the insight provided by a `coarse-grained' dynamics may be
sufficient for one's needs, specially when the actual, `sharp' dynamics is
too complex for a detailed analysis. In an \textit{ordinal representation},
a special case of symbolic dynamics, the state space is coarse-grained by
the ordinal partition $\pi _{L}$ (\ref{pi_L}), so the symbols are ordinal $L$%
-patterns. Some of these symbolic representations have an interesting
feature in common, namely, the alphabet (whether finite, countably infinite,
or continuous) is a group $\mathcal{G}$. Thus, $\mathcal{G}=(\mathbb{Z},+)$
when using bins, $\mathcal{G}=(\mathcal{S}_{L},\circ )$ (see (\ref{ros}))
when using ordinal $L$-patterns, and $\mathcal{G}=([0,1),+)$ when using
phases.

To show that group-theoretic representations can provide an additional
leverage in time series analysis, we consider next random processes with 
\textit{algebraic alphabets}, i.e., $\mathcal{G}$-valued random processes,
where $\mathcal{G}$ is a finite group.

\begin{definition}
The transfer entropy of two processes with a common algebraic alphabet is
called \emph{algebraic transfer entropy}.
\end{definition}

Actually, to warrant the convergence of information-theoretic quantities in
case of representations with infinite algebraic groups, it would be
sufficient that the probability of the symbols is positive only for a finite
number of them, as always happens in practice.

One way of exploiting the algebraic structure of the alphabet $\mathcal{G}$
is the following. Given $\alpha,\beta\in\mathcal{G}$ there exists a unique $%
\tau=\tau_{\alpha,\beta}\in\mathcal{G}$, called \textit{transcript }from the 
\textit{source symbol} $\alpha$ to the \textit{target symbol} $\beta$, such
that 
\begin{equation}
\tau=\beta\alpha^{-1},  \label{2_transcription}
\end{equation}
where, as usual, we denote the product of two elements by concatenating
them. The map $\mathcal{G}\times\mathcal{G}\rightarrow\mathcal{G}$ defined
by $(\alpha,\beta)\mapsto\tau_{\alpha,\beta}$ is $\left\vert \mathcal{G}%
\right\vert $-to-$1$ since the $\left\vert \mathcal{G}\right\vert $ distinct
pairs $(\alpha,\tau\alpha)$ are sent to $\tau$ for all $\alpha\in$ $\mathcal{%
G}$. Reciprocally, any pair $(\alpha,\beta)\in\mathcal{G\times G}$ whose
transcript is $\tau$ must have $\beta=\tau\alpha$ by (\ref{2_transcription}).

The \textit{coupling complexity coefficient} of the $\mathcal{G}$-valued
random variables $\alpha_{1}$,..., $\alpha_{N}$, denoted by $C(\alpha
_{1},...,\alpha_{N})$, is directly linked to their transcripts. They are
defined as \cite{Amigo2012A,Monetti2013A}%
\begin{equation}
C(\alpha_{1},\alpha_{2},...,\alpha_{N})=\min_{1\leq n\leq
N}I(\alpha_{n};\tau_{\alpha_{1},\alpha_{2}},\tau_{\alpha_{2},%
\alpha_{3}},...,\tau _{\alpha_{N-1},\alpha_{N}}).  \label{C()}
\end{equation}
By definition $C(\alpha_{1},...,\alpha_{N})\geq0$. The coefficients $%
C(\alpha_{1},...,\alpha_{N})$ have a number of interesting properties. A
basic one is the invariance under permutation of their arguments. See \cite%
{Monetti2013A,Amigo2014A} for other properties.

Let $\mathbf{\Xi}=(\xi_{t})_{t\in\mathbb{T}}$ and $\mathbf{\Upsilon}=(\eta
_{t})_{t\in\mathbb{T}}$ be two stationary, $\mathcal{G}$-valued random
processes (possibly obtained via group-theoretic representations of two
random processes $\mathbf{X}$ and $\mathbf{Y}$, respectively). The transfer
entropy from $\mathbf{\Upsilon}$\emph{\ }to\emph{\ }$\mathbf{\Xi}$ with
coupling delay $\Lambda\geq1$ is then (see (\ref{Transfer3})) 
\begin{equation}
T_{\mathbf{\Upsilon}\rightarrow\mathbf{\Xi}}(\Lambda)=I(\xi_{t+\Lambda};%
\eta_{t}\left\vert \xi_{t}\right) .  \label{transfer}
\end{equation}

\begin{theorem}
\cite{Amigo2014A} If $C(\xi_{t+\Lambda},\eta_{t},\xi_{t})=0$, and $H(\xi
_{t})\leq H(\eta_{t})$, then%
\begin{equation}
T_{\mathbf{\Upsilon}\rightarrow\mathbf{\Xi}}(\Lambda)=I(\tau_{\xi_{t+\Lambda
},\xi_{t}};\tau_{\eta_{t},\xi_{t}}).  \label{3.5b}
\end{equation}
\end{theorem}

Note that (\ref{3.5b}) equates a conditional mutual information (with three
variables and $\left\vert \mathcal{G}\right\vert ^{3}$ possible values) to
an unconditioned mutual information (with two variables and $\left\vert 
\mathcal{G}\right\vert ^{2}$ possible values) thanks to the use of
transcripts. This result, called the \textit{dimensional reduction of the
algebraic transfer entropy}, can make a difference in time series analysis
if the data sequences are short, as often happens in practice. See \cite[%
Theorem 1]{Amigo2014A} for much general results when the history lengths of
the processes involved are arbitrary.

\begin{corollary}
\label{Corollary1}If $C(\xi_{t+\Lambda},\eta_{t},\xi_{t})=C(\eta_{t+\Lambda
},\xi_{t},\eta_{t})=0$ and $H(\xi_{t})=H(\eta_{t})$, then one can calculate
the coupling direction indicator (\ref{DeltaTE}) as 
\begin{equation*}
\Delta T_{\mathbf{\Upsilon}\rightarrow\mathbf{\Xi}}(\Lambda)=I(\tau
_{\xi_{t+\Lambda},\xi_{t}};\tau_{\eta_{t},\xi_{t}})-I(\tau_{\eta_{t+\Lambda
},\eta_{t}};\tau_{\xi_{t},\eta_{t}})=-\Delta T_{\mathbf{\Xi}\rightarrow 
\mathbf{\Upsilon}}(\Lambda).
\end{equation*}
\end{corollary}

We conclude that the coupling direction between processes with a common
algebraic alphabet can be determined with mutual informations of
transcripts. Corollary \ref{Corollary1} was applied in \cite%
{Monetti2013B,Amigo2014A} to ordinal representations of time series analysis
with satisfactory results.

Transfer entropy remains a hot research topic because of the crucial role of
causality in data analysis. For further insights into this subtle subject,
the interested reader is referred to the papers \cite%
{Hahs2011,Runge2012,Smirnov2013,Liang2014,Cafaro2015}, just to mention a few.

\section{Applications}

The practical applications of the mathematical entropy are usually
associated with information theory and communications technology. Our pick
in this section belongs to the realm of applied mathematics. There are
certainly others in physics, computer science, and social sciences.

\subsection{Time series analysis}

Entropy plays a prominent role in the analy\-sis of time series, in
particular when the data stem from systems assumed to be non-linear. In this
context, a central objective of the data analyst is to quantify the
complexity of both data and systems. In the following we review some
applications of entropy to time series analysis. Needless to say, no claim
of completeness is made. Rather, we are going to concentrate on sample,
approximate, permutation, Tsallis and R\'{e}nyi entropies applied to
physiological data. Due to their very nature, these data contain a rich
structure of complex patterns, besides being often recorded in large
amounts. For applications of those entropies to the analysis of data from
other fields such as physics, chemistry, biology, economics and geophysics,
we refer to \cite%
{Tsallis2011,Gradojevic2011,Sneddon2007,BeadleEtAl2008,Amigo2015,AmigoKeller2013,ZaninZuninoRosso2012}
and the references therein.

Early applications of approximate, sample and permutation entropies to
phy\-siological data were directly connected with the formulation of the
corresponding concept or followed it a short time after (compare also \cite%
{RichmanMoorman2000,Pincus2006}). Thus, Pincus and Viscarello analyzed fetal
heart-rate variability (HRV) with approximate entropy already in 1992 \cite%
{Pincus1992}, and Lake et al., neonatal HRV with sample entropy in 2002 \cite%
{LakeRichmanGriffin2002}. Note that an abnormal heart rate is often
signalized by a reduced variability of the electrocardiograms (ECGs)
characteristics ---the right task for entropy! Frank et al. \cite%
{FrankPompeSchneider2006} demonstrated also the ability of permutation
entropy to classify fetal behavioral states on the basis of HRV. Many other
applications can be found in cardiovascular studies based on ECGs. For
instance, Acharya et al. \cite{AcharyaEtAl2006} and Voss et al. \cite%
{VossSchulz2009} compared the performances of various measures, including
approximate and sample entropy, when analyzing HRV. Along the same line,
Graff et al. \cite{GraffGraffKaczkowska2012} studied the ability of the
approximate, sample, permutation, and fuzzy entropies to discriminate
between healthy patients and patients with congestive heart failure in the
case of short ECG data sets. For further applications and comparisons of
entropies related to ECGs, we refer to \cite{Amigo2015,ZaninZuninoRosso2012}.

Just as ECGs are indispensable for measuring heart activity,
electroencephalograms (EEGs) and magnetoencephalograms (MEGs) are the
principal diagnostic tools for assessing brain activity. Since abnormal
brain states are often reflected by a change in the complexity of the
electromagnetic activity, the interest in entropy for the analysis of EEGs
and MEGs is increasing in the medical research and praxis.

Approximate, sample and permutation entropies have been the basic tools in
the non-linear analysis of EEGs and MEGs from patients with Alzheimer's
disease (AD), mainly to distinguish healthy from diseased subjects (see
Hornero et al. \cite{HorneroAbasoloEscuredo2009}, Morabito et al. \cite%
{Morabito2012}, and references therein). Another field of applications is
the quantification of anesthetic drug effects on the brain activity as
measured by EEGs, including comparative testing of different anesthetics and
the discrimination between consciousness and unconsciousness. Thereby
different entropy measures have been shown to be sensitive to the depth of
anaesthesia. For an overview on the main entropy measures used in this field
and a comprehensive list of references, see Liang et al. \cite{LiangEtAl2015}%
.

Some additional applications of the approximate and sample entropies are
related to sleep research, in particular to the separation of sleep stages
based on EEG data. For example, Acharya et al. \cite{AcharyaEtAl2005}
compared several non-linear measures, including approximate entropy, for
analyzing sleep stages in surface EEGs, while Burioka et al. \cite%
{BuriokaMiyataCornelissen2005} estimated the values of the approximate
entropy for different sleep stages. Also permutation entropy has been
applied to classify sleep stages in human EEGs; the first attempt in this
direction was made in \cite{NicolaouGeorgiou2011}. Moreover, sleep stage
segmentation via the conditional entropy of ordinal patterns (directly
related to permutation entropy) has been studied in \cite%
{KellerUnakafovUnakafova2014}.

The epilepsy research is probably the main biomedical application field of
entropy in time series analysis. The reasons for this are manifold. About
one percent of the world's population is estimated to suffer from epilepsy
in a spectrum that goes from mild to strong forms, the latter case being
associated with serious health and psychosocial problems. Moreover, the
mechanisms of epileptic activity are far from being understood. From the
viewpoint of non-linear time series analysis, EEG signals related to
epileptic activity are interesting because they contain special wave forms
(spike waves, sharp waves, slow waves), often believed to be the result of
low-dimensional dynamics.

The approximate and sample entropies have been used as biomarkers in
algorithms for epileptic EEG analysis, in particular for epileptic seizure
detection. In this category, Kannatal et al. \cite{Kannathal2005} tested
different entropy measures, among them approximate entropy, using a popular
data set related to epilepsy \cite{AndrzejakLehnertzMormann2001}. Srinivasan
et al. \cite{SrinivasanEswaranSriraam2007} investigated approximate entropy
as an input feature for a neural-network-based automated epileptic EEG
detection system, and Jouny et al. \cite{JounyBergey2012} considered
different complexity and entropy measures, in particular sample entropy and
permutation entropy, to characterize early partial seizure onset in
intracranial EEG recordings.

The first application of permutation entropy in epileptic EEG analysis is
due to Cao et al. \cite{CaoTungGao2004}, who detected this way dynamic
changes during epileptic seizures in intracranial EEGs, and also to Keller
and Lauffer \cite{KellerLaufer2003}, who studied vagus stimulation for
reducing epileptic activity. The predictability of epileptic seizures by
permutation entropy has been the subject of several papers (see, e.g., Li et
al. \cite{LiOuyangRichards2007} and Bruzzo et al.\cite%
{BruzzoGesierichSanti2008}). Li et al. \cite{LiYanLiu2014} studied how
permutation entropy reflects the changes in surface EEGs due to absence
seizures during different seizure phases (seizure-free, pre-seizure,
seizure). The choice of parameters in applications of permutation entropy to
time series from different sleep stages were studied in \cite%
{KellerUnakafovUnakafova2014,UnakafovaThesis}.

As discussed in Sect. 2, the Shannon entropy along with its generalizations,
the Tsallis and R\'{e}nyi entropies, are functions of discrete probability
distributions, which for time series are not given a priori. The $q$%
-dependent weighting of the probabilities in the Tsallis and R\'{e}nyi
entropies as compared to the Shannon entropies endows data analysis with
more flexibility. There are various ways of extracting probability
distributions from time series. One of them (already mentioned in previous
sections) is symbolic dynamics, which in practice boils down to counting
symbols or patterns in time series. This being the case, it is natural to
consider the Tsallis and R\'{e}nyi variants of permutation entropy. First
attempts in this direction have been reported in \cite{LiangEtAl2015} in the
field of anesthesia. Moreover, R\'{e}nyi entropy in combination with simple
symbolic dynamics has been applied to HRV by Kurths et al.; see, e.g., \cite%
{Kurths1995}.

Furthermore, the estimation of the probability distributions for computing
entropies can be done in the frequency domain as well. The interest shifts
now to the distributions of frequencies and wavelets at different scales. In
this context, the Tsallis entropy has been used for EEG and ECG analysis by
a number of groups, e.g., Gamero et al. \cite{GameroPlastinoTorres1997},
Capurro et al. \cite{CapurroDiambraLorenzo1998}, and Zhang et al. \cite%
{ZhangEtAl2010}). As for the R\'{e}nyi entropy, applications reach from
epilepsy detection in EEG (see, e.g., \cite{Kannathal2005}) over artifact
rejection in multichannel scalp EEG (see \cite{Mammone2012} and references
therein) to early diagnosis of Alzheimer's disease in MEG data (see, e.g., 
\cite{Poza2008}).

\subsection{Statistical inference}

Let $X$ be a random variable on a probability space $(\Omega ,\mathcal{B}%
,\mu )$ taking values in a measurable set $\Gamma \subset \mathbb{R}$.
Suppose that all we know about $X$ are the expectation values of some
measurable functions (`observables') $\phi _{k}:\Gamma \rightarrow \mathbb{R}
$, $1\leq k\leq K$, i.e.,%
\begin{equation}
\left\langle \phi _{k}(X)\right\rangle :=\int_{S}\phi _{k}(x)\rho (x)\mathrm{%
d}x=m_{k},  \label{moments}
\end{equation}%
where $\rho (x)\geq 0$ is the probability density function of $X$, and $%
S\subset \Gamma $ is the support of $\rho (x)$ (i.e., $S=\{x\in \Gamma :\rho
(x)>0\}$), hence%
\begin{equation}
\int_{S}\rho (x)\mathrm{d}x=1.  \label{rho}
\end{equation}%
The conditions (\ref{moments}) are called the \textit{moment constraints}.
To guarantee that the integrals (\ref{moments}) exist, we suppose that the
functions $\phi _{k}$, $1\leq k\leq K$, are integrable over $S$.

Given the functions $\phi _{k}$ and the numbers $m_{k}$, $0\leq k\leq K$,
how can we use this information to best characterize the density function $%
\rho (x)$ of $X$?

According to the \textit{maximum entropy principle} of Jaynes \cite%
{Jaynes1957}, \textquotedblleft in making inferences on the basis of partial
information we must use the probability distribution which has maximum
entropy subject to whatever is known\textquotedblright . In the case that
nothing is known, this principle leads to Laplace's principle of
indifference, also called the \textit{principle of insufficient reason},
according to which all events will be assigned equal probabilities.
Therefore, the principle of maximum entropy can be considered as an
extension of Laplace's principle.

\begin{theorem}
\label{Thm_maxent}\emph{(Maximum entropy probability distributions)} The
probability density function $\rho ^{\ast }(x)$ which maximizes the entropy%
\begin{equation}
H(\rho )=-\int_{S}\rho (x)\ln \rho (x)\mathrm{d}x  \label{H(rho)}
\end{equation}%
over all probability density functions satisfying the moment constraints (%
\ref{moments}) are of the form%
\begin{equation*}
\rho ^{\ast }(x)=\exp \left( \lambda _{0}-1+\sum_{k=1}^{K}\lambda _{k}\phi
_{k}(x)\right) ,
\end{equation*}%
$x\in S$, provided that there exist $\lambda _{0},...,\lambda _{K}\in 
\mathbb{R}$ such that $\rho ^{\ast }(x)$ satisfies the constraints (\ref%
{moments}) and (\ref{rho}).
\end{theorem}

See \cite[Chapter 12]{Cover2006} for a proof without resorting to calculus
of variations. Needless to say, the parameters $\lambda_{j}$ are the
Lagrange multipliers corresponding to the constraints (\ref{rho}) for $j=0$,
and (\ref{moments}) for $1\leq j\leq K$.

In \cite{Park2009} one can find a table with a number of maximum entropy
probability distributions corresponding to a variety of moment constraints.
We summarize next the simplest cases: with no moment constraints, with fixed
mean, and with fixed mean and variance.

\begin{corollary}
\label{CorollaryG}\emph{(Basic maximum entropy distributions)}

\begin{enumerate}
\item If $S=[a,b]$, and there are no moment constraints, then $%
\rho^{\ast}(x) $ is the uniform distribution over $[a,b]$, i.e., $%
\rho^{\ast}(x)=\frac{1}{b-a}$, and $H(\rho^{\ast})=-\ln(b-a)$.

\item If $S=[0,+\infty)$ and $\left\langle X\right\rangle =m$, then $%
\rho^{\ast}(x)$ is the exponential distribution%
\begin{equation*}
\rho^{\ast}(x)=\frac{1}{m}e^{-\frac{x}{m}},
\end{equation*}
and $H(\rho^{\ast})=\ln m+1$. Otherwise, if no moment constraint is given,
then there is no maximum entropy probability distribution.

\item If $S=(-\infty,+\infty)$ and no moment constraints are given, or only
the mean is prescribed, then there is no maximum entropy probability
density. If the moments $\left\langle X\right\rangle =m_{1}$ and $%
\left\langle X^{2}\right\rangle =m_{2}$ are given, then the maximum entropy
probability distribution is the normal distribution $\mathcal{N}%
(m_{1},m_{2}-m_{1}^{2})$.
\end{enumerate}
\end{corollary}

At variance with Corollary \ref{CorollaryG}(3), maximization of the entropy
with prescribed $\left\langle X\right\rangle =m_{1}$, $\left\langle
X^{2}\right\rangle =m_{2}$, and $\left\langle X^{3}\right\rangle =m_{3}$ in $%
S=(-\infty ,+\infty )$ is, in general, not possible \cite[Chapter 12]%
{Cover2006}, although one can come arbitrarily close the upper bound for the
maximum entropy distribution 
\begin{equation*}
\sup H(\rho )=H(\mathcal{N}(0,m_{2}-m_{1}^{2}))=\tfrac{1}{2}\ln 2\pi
e(m_{2}-m_{1}^{2}),
\end{equation*}%
More generally, it can happen that $\sup H(\rho )$ is not achievable with
probability densities\ satisfying all the constraints.

The same analysis holds for multivariate random variables. Thus, if $S=%
\mathbb{R}^{n}$ and $\left\langle X_{i}X_{j}\right\rangle =R_{ij}$, $1\leq
i,j\leq n$, then the maximum entropy density is the (zero mean) Gaussian
distribution%
\begin{equation*}
\rho^{\ast}(\mathbf{x})=\frac{1}{\sqrt{(2\pi)^{n}\det R}}\exp\left( -\tfrac{1%
}{2}\mathbf{x}^{T}R^{-1}\mathbf{x}\right) ,
\end{equation*}
where $\mathbf{x}=(x_{1},...,x_{n})^{T}$ and $R$ is the autocorrelation
matrix $(R_{ij})_{1\leq i,j\leq n}$. One finds%
\begin{equation*}
H(\rho^{\ast})=\frac{1}{2}\ln\left( (2\pi e)^{n}\det R\right) .
\end{equation*}

Random variables with finite alphabets are a particular case of Theorem \ref%
{Thm_maxent}. If $\Gamma =\{x_{1},...,x_{N}\}$, and $\mu (X=x_{n})=p_{n}$, $%
1\leq n\leq N$, then the maximum entropy probability distribution is%
\begin{equation}
p_{n}^{\ast }=\frac{1}{Z(\lambda _{1},...,\lambda _{K})}\exp \left(
\sum_{i=1}^{K}\lambda _{i}\phi _{i}(x_{n})\right) ,  \label{p_i}
\end{equation}%
where the normalization factor 
\begin{equation}
Z(\lambda _{1},...,\lambda _{K})=e^{1-\lambda _{0}}=\sum_{n=1}^{N}\exp
\left( \sum_{i=1}^{K}\lambda _{i}\phi _{i}(x_{n})\right)  \label{Z()}
\end{equation}%
is called the \textit{partition function}. It follows%
\begin{align}
m_{k}& =\sum_{n=1}^{N}p_{n}^{\ast }\phi _{k}(x_{n})=\frac{1}{Z(\lambda
_{1},...,\lambda _{K})}\sum_{n=1}^{N}\phi _{k}(x_{n})\exp \left(
\sum_{i=1}^{K}\lambda _{i}\phi _{i}(x_{n})\right)  \label{e_k} \\
& =\frac{\partial }{\partial \lambda _{k}}\ln Z(\lambda _{1},...,\lambda
_{K})  \notag
\end{align}%
for $1\leq k\leq K$. Equations (\ref{e_k}) determine the Lagrange
multipliers $\lambda _{1},...,\lambda _{K}$, provided that these can be
solved as functions of $m_{1},...,m_{K}$.

The probability distribution that maximizes the Gibbs entropy (\ref{Gibbs_S}%
) subject to the energy constraint (\ref{canonical}) is the so-called Gibbs
distribution for the canonical ensemble,%
\begin{equation*}
p_{i}^{\ast }\sim \exp \left( -\frac{1}{k_{B}T}\varepsilon _{i}\right) .
\end{equation*}%
Here $\Gamma =\{\varepsilon _{1},...,\varepsilon _{N}\}$, $\phi _{1}(x)=x$, $%
m_{1}=U$, and $\lambda _{1}=-\frac{1}{k_{B}T}$. If instead of the Gibbs
entropy one uses the Tsallis entropy $S_{q}$, Equation (\ref{Tsallis}), then 
\begin{equation*}
p_{i}^{\ast }(q)\sim \left( 1-\frac{q-1}{k_{B}T}\varepsilon _{i}\right) ^{%
\frac{1}{q-1}}.
\end{equation*}

The fact that the probability distributions of the statistical mechanics
(not only for the canonical ensemble) can be deduced by
information-theoretical means led Jaynes to conclude that statistical
mechanics is an application of information theory. More recent and
challenging applications of entropy maximization to collective behavior
(e.g., tool use and social cooperation) can be found in \cite{Wissner2013}
and the references therein.

\subsection{Cluster analysis}

Application of entropy with its interpretation as amount of information,
uncertainty or diversity has led also to various kinds of minimum entropy
principles in different fields. As representatives of such principles, we
are going to touch upon two strategies used in data clustering.

In cluster analysis one seeks a natural partitioning of a set of $n$ objects
into subsets $c_{1},c_{2},\ldots ,c_{k}$ called the \textit{clusters}. For
simplicity, we assume the objects to be identified with the \textit{feature
vectors} $x_{1},x_{2},\ldots ,x_{n}\in {\mathbb{R}}^{m}$, whose components
are the values of $m$ observables at each object. Clearly, the question what
a natural partition is, strongly depends on the practical aspects of a data
analysis, with the result that the number of different clustering criteria
is huge. Nevertheless, a rough criterion is that the feature vectors inside
a cluster should be similar, while dissimilar when belonging to different
clusters. An information-theoretical tool to describe the homogeneity of
data within a cluster is, of course, entropy. In this context, the idea of
minimizing the mean entropy of clusters is straightforward.

Assume that the feature vectors are obtained by independent draws of a
continuous $m$-dimensional random vector $X$, and describe the cluster
assignment by a random variable $C$ taking on the values $i=1,2,\ldots ,n$,
the outcome $i$ being the label of the cluster a feature vector belongs to.
This way, a realization of the model $(X,C)$ can be identified with a
feature vector and the cluster containing it, respectively. In this setting,
the minimalization of the mean cluster entropy amounts to searching for a
cluster assignment $C$ with minimal conditional entropy 
\begin{eqnarray*}
H(X|C)=-\sum_{i=1}^n p(i)\int \rho(x|i)\ln \rho(x|i)\,dx,
\end{eqnarray*}
where $p(i)$ denotes the probability of the cluster $c_i$ and $\rho(x|i)$
the conditional density of $X$ given $C=i$ (see \cite{RobertsEtAl1999}). In
a sense, one seeks a cluster assignment $C$ for which the variable $X$
representing the feature vectors contains as much information of $C$ as
possible.

From an automatic learning viewpoint, it is more natural to start from the
feature vectors and look for a clustering optimally representing the
structure of the set of feature vectors. This approach, called the \textit{%
minimum conditional entropy principle} (see e.g.~\cite{LiEtAl2004,DaiHu2010}%
), consists in minimizing the conditional entropy 
\begin{eqnarray*}
H(C|X)=-\int \sum_{i=1}^n \rho(i|x)\ln \rho(i|x)\rho(x)\,d\rho(x),
\end{eqnarray*}
where $\rho(x)\geq0$ is the probability density function of $X$ and $%
\rho(i|x)$ is the conditional probability of $i$ given $x$.

Both models, possibly together with additional assumptions about the
probability distributions, lead to clustering strategies for given feature
vectors $x_{1},x_{2},\ldots ,x_{n}$ via minimizing appropriate estimations
of $H(X|C)$ or $H(C|X)$ on the base of the feature vectors. Note that in 
\cite{LiEtAl2004} a Tsallis variant of $H(C|X)$ has been also discussed.

\subsection{Testing stochastic independence and determinism}

To finish this sample of applications, let us discuss some stochastic
independence tests based both on metric and topological permutation entropy.
Along with the biomedical applications to time series analysis mentioned in
Sect. 5.1, testing financial time series for stochastic independence was
also one of the first applications of permutation entropy \cite%
{Matilla2007,Matilla2008}.

As already explained in Sect. 4.3, an immediate consequence of Theorem 4.11%
\textit{(ii)} for the topological permutation entropy is the existence of
forbidden ordinal patterns in the restricted setting of Definition 4.10. To
be specific, if $I$ is a closed interval of $\mathbb{R}$ and $T:I\rightarrow
I$ is a piecewise, strictly monotone map, then there exist $L_{\min }\geq 2$
and length-$L$ ordinal patterns $r$ such that no $x\in I$ is of type $r$ for
all $L\geq L_{\min }$. In other words, with the mild provisos just stated
(which may be taken for granted in practice), a one-dimensional dynamics has
always forbidden ordinal patterns of sufficiently large length. The opposite
happens with independent and identically distributed (i.i.d.) random
processes, also called \textit{white noise}. In this case, all ordinal
patterns of whatever length are allowed.

Let $(x_{t})_{t\in \mathbb{N}_{0}}$ be the orbit of $x_{0}\in I$ under the
action of $T$. In practice, an observer will not measure the true values $%
x_{t}$ but the `noisy' values $\tilde{x}_{t}$ instead, due to a number of
reasons such as measurement errors and electronic noise. This fact is
modeled by adding to $x_{t}$ a so-called observational noise $w_{t}$ at each
time step, i.e.,%
\begin{equation}
\tilde{x}_{t}=x_{t}+w_{t},  \label{noisy}
\end{equation}%
where $(w_{t})_{t\in \mathbb{N}_{0}}$ is a realization of a stochastic
process. Unless one knows better, the $w_{t}$'s are supposed to be white
noise. If, on the contrary, there are correlations among the variables $%
w_{t} $, then one speaks of colored noise.

One of the preliminary tasks in non-linear time series analysis is precisely
to ascertain whether the noisy data at hand has actually been output by a
deterministic system, i.e., if they have the structure (\ref{noisy}). There
are some standard methods for detecting determinism, e.g., the cross
prediction error statistic \cite{Kantz2005}, but denoising should be first
applied.

An alternative method in one-dimensional dynamics, which does not require
denoising, exploits the existence of forbidden ordinal patterns, and the
robustness of ordinal patterns against observational noise \cite%
{Amigo2008B,Amigo2010A,Amigo2010B}. To this end, we use hypothesis testing.
Our \textit{null hypothesis} is that the noisy data $(\tilde{x}_{t})_{t\in 
\mathbb{N}_{0}}$ are outcomes of an i.i.d random process, that is, 
\begin{equation}
H_{0}:\text{ }(\tilde{x}_{t})_{t\in \mathbb{N}_{0}}\text{ is white noise.}
\label{H_0}
\end{equation}%
A working hypothesis in non-linear time series analysis (based on the
phenomenon being observed) is that, however random the data may look like,
there is an underlying deterministic component. This being the case, the
rejection of $H_{0}$ is equated to determinism.

A further caveat is that, in practice, all time series are finite. This
implies the possible existence of false forbidden patterns, i.e., allowed
ordinal patterns that are missing just because the time series is finite but
would appear if the series were long enough \cite{Amigo2007}. Therefore, it
is always good practice to check the stability of the results against
changes in the length of the data.

Thus, consider a noisy and finite time series $(\tilde{x}_{t})_{t=0}^{N-1}$
and the null hypothesis%
\begin{equation}
\tilde{H}_{0}:\text{ }(\tilde{x}_{t})_{t=0}^{N-1}\text{ is i.i.d.}
\label{H_0B}
\end{equation}
To accept or reject $\tilde{H}_{0}$, there are two simple-minded methods.

\textbf{Method 1} \cite[Chapter 9]{Amigo2010A} consists of the following
three steps.

\begin{description}
\item[(a)] Count the number of missing ordinal $L$-patterns of adequate
length (say $(L+1)!\leq N$) using a sliding window of size $L$.

\item[(b)] Randomize the time series, i.e., create a sample of surrogate
data.

\item[(c)] Proceed as in step (a) with the surrogates.
\end{description}

Use now any measure of dissimilarity for the number of forbidden patterns
(or even for the probability distributions) to accept or reject $\tilde{H}%
_{0}$. Roughly speaking, if the results of (a) and (c) are about the same,
the sequence is very likely not deterministic (or the observational noise is
so strong as compared to the deterministic signal that the latter has been
completely masked). Otherwise, the data stems from a deterministic process.

\textbf{Method 2} \cite[Chapter 9]{Amigo2010A} uses also sliding windows,
but this time without overlap. The number of such windows is $K=\left\lfloor
N/L\right\rfloor $, each comprising the entries%
\begin{equation*}
\mathbf{e}_{k}=(\tilde{x}_{kL},\tilde{x}_{kL+1},...,\tilde{x}%
_{(k+1)L-1}),\;\;0\leq k\leq K-1.
\end{equation*}%
Under the null hypothesis, the ordinal patterns of the non-overlapping
windows $\mathbf{e}_{k}$ will be also i.i.d.. Let $\nu _{r}$ be the number
of $\mathbf{e}_{k}$'s of type $r\in \mathcal{S}_{L}$. Thus, $\nu _{r}=0$
means that the $L$-pattern $r$ is missing in $(\tilde{x}_{t})_{t=0}^{N-1}$;
otherwise we say that $r$ is visible. We apply then a chi-square
goodness-of-fit hypothesis test with the statistic%
\begin{equation*}
\chi ^{2}(L)=\sum_{r\in \mathcal{S}_{L}}\frac{(\nu _{r}-K/L!)^{2}}{K/L!}=%
\frac{L!}{K}\sum_{r\in \mathcal{S}_{L}:\text{visible}}\nu _{r}^{2}-K,
\end{equation*}%
since $\sum_{\pi \in \mathcal{S}_{L}}\nu _{\pi }=K$. Here $K/L!$ is the
absolute frequency of an ordinal $L$-pattern if $\tilde{H}_{0}$ holds true.
In the affirmative case, $\chi ^{2}(L)$ converges in distribution (as $%
K\rightarrow \infty $) to a chi-square distribution with $L!-1$ degrees of
freedom. Thus, for large $K$, a test with approximate level $\alpha $ is
obtained by rejecting $\tilde{H}_{0}$ if $\chi ^{2}(L)>\chi _{L!-1,1-\alpha
}^{2}$, where $\chi _{L!-1,1-\alpha }^{2}$ is the upper $1-\alpha $ critical
point for the chi-square distribution with $L!-1$ degrees of freedom \cite%
{Law2000}. In our case, the hypothetical convergence of $\chi ^{2}(L)$ to
the corresponding chi-square distribution may be considered sufficiently
good if $\nu _{\pi }>10$ for all visible $L$-patterns $\pi $, and $K/L!>5$.

Notice that, since this test is based on distributions, it would be possible
that a deterministic map has no forbidden $L$-patterns, thus $\nu _{\pi
}\neq 0$ for all $\pi \in \mathcal{S}_{L}$, however, the null hypothesis
could be rejected because those $\nu _{\pi }$'s are not evenly distributed.

A standard test for stochastic independence in time series is the
Brock-Dechert-Scheinkman (DBS) test, which is based on the correlation
dimension. Method 2 was numerically benchmarked against the DBS test in \cite%
{Amigo2008B,Amigo2010B} by means of the Lorenz map and the delayed Henon
map. In practically all cases studied, Method 2 outperformed the DBS test. A
similar method based on a different statistic was proposed in \cite%
{Canovas2009}.

Finally, one can also use the metric permutation entropy to test the null
hypothesis (\ref{H_0B}). To this end, let%
\begin{equation*}
h_{L}^{\ast }((\tilde{x}_{t})_{t=0}^{N-1})=-\sum_{r\in \mathcal{S}%
_{L}}p_{r}\ln p_{r}
\end{equation*}%
be the empirical metric permutation entropy of order $L$ of $%
(x_{t})_{t=0}^{N-1}$, i.e.,%
\begin{equation*}
p_{r}=\frac{\left\vert \{(\tilde{x}_{t},\tilde{x}_{t+1},...,\tilde{x}%
_{t+L-1})\text{ of type }r\in \mathcal{S}_{L}:0\leq t\leq N-L\}\right\vert }{%
N-L+1}.
\end{equation*}

\medskip\ 

\noindent \textbf{Theorem} \cite{Matilla2008}. If $(\tilde{x}%
_{t})_{t=0}^{N-1}$ is i.i.d., then the statistic%
\begin{equation*}
G(L)=2(N-L+1)\left( \ln L!-h_{L}^{\ast }((\tilde{x}_{t})_{t=0}^{N-1})\right)
\end{equation*}%
is asymptotically $\chi _{L!-1}^{2}$ distributed.

\medskip

Therefore, the decision rule at a $100(1-\alpha )\%$ confidence level is the
following: If $G(L)\leq \chi _{L!-1,\alpha }^{2}$ accept $\tilde{H}_{0}$,
otherwise reject $\tilde{H}_{0}$.

\bigskip

\section{Acknowledgements}

This work was financially support by the Spanish \textit{Ministry of Economy
and Competitivity}, grant MTM2012-31698, and the \textit{Graduate School for
Computing in Medicine and Life Sciences} funded by \textit{Germany
Excellence Initiative} [DFG GSC 235/1].


\bibliographystyle{unsrt}
\bibliography{DCDSB2015AArxiv5}

%
%
%
%

\end{document}